\documentclass[12pt]{amsart}
\usepackage{amssymb,amsmath,amsfonts,latexsym}
\usepackage{bm}
\usepackage[all,cmtip]{xy}
\usepackage{amscd}
\usepackage{multirow,bigdelim}

\newcommand{\newsection}[1]{\setcounter{equation}{0} \section{#1}}
\newcommand{\bea}{\begin{eqnarray}}
\newcommand{\eea}{\end{eqnarray}}

\newcommand{\vp}{\varphi}

\newcommand{\cla}{\mathcal{A}}
\newcommand{\clb}{\mathcal{B}}

\newcommand{\cle}{\mathcal{E}}
\newcommand{\clf}{\mathcal{F}}

\newcommand{\clh}{\mathcal{H}}
\newcommand{\clk}{\mathcal{K}}
\newcommand{\cll}{\mathcal{L}}
\newcommand{\clm}{\mathcal{M}}
\newcommand{\cln}{\mathcal{N}}

\newcommand{\cls}{\mathcal{S}}

\newcommand{\z}{\bm{z}}
\newcommand{\w}{\bm{w}}

\newcommand{\D}{\mathbb{D}}
\newcommand{\ball}{\mathbb{B}}

\newcommand{\raro}{\rightarrow}

\def\textmatrix#1&#2\\#3&#4\\{\bigl({#1 \atop #3}\ {#2 \atop #4}\bigr)}
\def\dispmatrix#1&#2\\#3&#4\\{\left({#1 \atop #3}\ {#2 \atop #4}\right)}
\newcommand{\be}{\begin{equation}}
\newcommand{\ee}{\end{equation}}
\newcommand{\ben}{\begin{eqnarray*}}
\newcommand{\een}{\end{eqnarray*}}

\newcommand{\NI}{\noindent}

\newcommand{\bi}{\begin{itemize}}
\newcommand{\ei}{\end{itemize}}

\theoremstyle{definition}

\theoremstyle{plain}

\newtheorem{thm}{Theorem}[section]

\newtheorem{defn}[thm]{Definition}
\newtheorem{rem}[thm]{Remark}

\numberwithin{equation}{section}

\let\phi=\varphi

\begin{document}

\title[Factorizations of Schur functions]{Factorizations of Schur functions}

\author[Debnath]{Ramlal Debnath}
\address{Indian Statistical Institute, Statistics and Mathematics Unit, 8th Mile, Mysore Road, Bangalore, 560059, India}
\email{ramlal\_rs@isibang.ac.in, ramlaldebnath100@gmail.com}

\author[Sarkar]{Jaydeb Sarkar}
\address{Indian Statistical Institute, Statistics and Mathematics Unit, 8th Mile, Mysore Road, Bangalore, 560059, India}
\email{jay@isibang.ac.in, jaydeb@gmail.com}





\begin{abstract}
The Schur class, denoted by $\mathcal{S}(\mathbb{D})$, is the set of all functions analytic and bounded by one in modulus in the open unit disc $\mathbb{D}$ in the complex plane $\mathbb{C}$, that is
\[
\mathcal{S}(\mathbb{D}) = \{\varphi \in H^\infty(\mathbb{D}): \|\varphi\|_{\infty} := \sup_{z \in \mathbb{D}} |\varphi(z)| \leq 1\}.
\]
The elements of $\mathcal{S}(\mathbb{D})$ are called Schur functions. A classical result going back to I. Schur states: A function $\varphi: \mathbb{D} \rightarrow \mathbb{C}$ is in $\mathcal{S}(\mathbb{D})$ if and only if there exist a Hilbert space $\mathcal{H}$ and an isometry (known as colligation operator matrix or scattering operator matrix)
\[
V = \begin{bmatrix} a & B \\ C & D \end{bmatrix} : \mathbb{C} \oplus \mathcal{H} \rightarrow \mathbb{C} \oplus \mathcal{H},
\]
such that $\varphi$ admits a transfer function realization corresponding to $V$, that is
\[
\varphi(z) = a + z B (I_{\mathcal{H}} - z D)^{-1} C \quad \quad (z \in \mathbb{D}).
\]
An analogous statement holds true for Schur functions on the bidisc. On the other hand, Schur-Agler class functions on the unit polydisc in $\mathbb{C}^n$ is a well-known ``analogue'' of Schur functions on $\mathbb{D}$. In this paper, we present algorithms to factorize Schur functions and Schur-Agler class functions in terms of colligation matrices. More precisely, we isolate checkable conditions on colligation matrices that ensure the existence of Schur (Schur-Agler class) factors of a Schur (Schur-Agler class) function and vice versa.
\end{abstract}

\maketitle

\newsection{Introduction}

In this paper, $\D^n$ denotes the open unit polydisc in $\mathbb{C}^n$, $n \geq 1$. By definition, the classical \textit{Schur class} $\cls(\D^n)$ consists of complex-valued analytic functions mapping from $\D^n$ into the closed unit disk $\overline{\D}$, that is
\[
\cls(\D^n) = \{\varphi : \D^n \raro \mathbb{C}: \varphi \mbox{~is analytic and ~} \|\varphi\|_{\infty} \leq 1\},
\]
where $\|\cdot\|_{\infty}$ denotes the supremum norm over $\D^n$. In other words, $\cls(\D^n)$ is the closed unit ball of the commutative Banach algebra $H^\infty(\D^n)$, the set of all bounded analytic functions on $\D^n$ under the supremum norm. The elements in the set $\mathcal{S}(\mathbb{D}^n)$ are called Schur functions \cite{Schur1, Schur2}.

It is a very remarkable fact that the one variable (and two variables too but not more than two variables, as we will see soon) Schur functions are closely related, via isometric colligations (or ``lurking isometries'' \cite{AM}), to bounded linear operators on Hilbert spaces. Recall that a \textit{colligation} (or \textit{scattering operator matrix}) is any bounded linear operator $V$ of the form
\[
V = \begin{bmatrix} A & B \\ C & D \end{bmatrix} : \cle \oplus \clh \raro \cle_* \oplus \clh,
\]
where $\clh$, $\cle$ and $\cle_*$ are Hilbert spaces. The colligation is said to be isometry if $V$ is isometry. Now, let $\clh$ be a Hilbert space and let
\begin{equation}\label{eq-1 var colligation}
V = \begin{bmatrix} a & B \\ C & D \end{bmatrix} : \mathbb{C} \oplus \clh \raro \mathbb{C} \oplus \clh,
\end{equation}
be an isometric colligation. Then a straightforward but lengthy and conceptual calculation (cf. page 73, \cite{AM}) verifies that $\tau_V \in \cls(\D)$, where
\[
\tau_V(z) = a + z B (I_{\mathcal{H}} - z D)^{-1} C \quad \quad (z \in \D).
\]
We call $\tau_V$ the \textit{transfer function realization} of the isometric colligation $V$. Conversely, if $\vp \in \cls(\D)$, then there exist a Hilbert space $\clh$ and an isometric colligation $V$ on $\mathbb{C} \oplus \clh$, as in \eqref{eq-1 var colligation}, such that
\[
\vp = \tau_V.
\]

We now pause, with our background so far, to state one of our main results specializing to the $n = 1$ case (see Theorem \ref{th-factor 3}): Suppose $\vp \in \cls(\D)$. If $\vp = \vp_1 \vp_2$ for some $\vp_1$ and $\vp_2$ in $\cls(\D)$, then there exist Hilbert spaces $\clh_1$ and $\clh_2$ and an (explicit) isometric colligation
\begin{equation}\label{eq-sec ex1}
V=\begin{bmatrix}a&B\\C&D\end{bmatrix} :=
\left[\begin{array}{@{}c|cc@{}} a&B_1&B_2\\\hline
C_1& D_{11}&D_{12}\\
C_2&D_{21}&D_{22}
\end{array}\right] :\mathbb{C} \oplus (\mathcal{H}_1\oplus\mathcal{H}_2) \rightarrow \mathbb{C} \oplus (\mathcal{H}_1\oplus\mathcal{H}_2),
\end{equation}
such that
\begin{equation}\label{eq-sec ex2}
D_{21}=0 \quad \mbox{and} \quad aD_{12} = C_1B_2,
\end{equation}
and $\vp = \tau_V,$
where $\tau_V(z) = a + z B (I_{\clh_1 \oplus \clh_2} - z D)^{-1} C$, $z \in \D$.
\\
The converse is true under an additional assumption that $\vp(0) \neq 0$ (see Theorem \ref{thm-new 5 case1}, Section \ref{sbusect theta 0}, for the case $\vp(0) = 0$): If $\vp = \tau_V$ for some isometric colligation $V$ as in \eqref{eq-sec ex1} satisfying \eqref{eq-sec ex2} and $a :=\vp(0) \neq 0$, then $\vp = \vp_1 \vp_2$ for some $\vp_1$ and $\vp_2$ in $\cls(\D)$. Moreover, in this case, $\phi$ and $\psi$ are explicitly given by $\vp_1 = \tau_{V_1}$ and $\vp_2 = \tau_{V_2}$ where
\[
V_1 = \begin{bmatrix}\alpha & B_1 \\ \frac{1}{\beta} C_1 & D_{11} \end{bmatrix} \in \clb(\mathbb{C} \oplus \clh_1) \quad \mbox{and} \quad V_2 = \begin{bmatrix}\beta & \frac{1}{\alpha} B_2 \\ C_2 & D_{22} \end{bmatrix} \in \clb(\mathbb{C} \oplus \clh_2),
\]
are isometric colligations and $\alpha$ and $\beta$ are non-zero scalars which satisfy the following conditions
\[
|\beta|^2 = |a|^2 + C_1^* C_1 \quad \mbox{and} \quad \alpha = \frac{a}{\beta}.
\]

In view of the above results, it is now clear that the goal of this paper is to clarify the link between isometric colligations and factors of Schur functions.

We also remark that the above one-variable factorization of Schur functions also relates to factorizations of Sz.-Nagy and Foias characteristic functions \cite{Helton} as well as Brodski\u{i} colligations \cite{Brod} in terms of invariant subspaces of certain operators \cite[Theorem 2.6]{Brod}. More specifically, see the idea of the product of colligations (as well as for a similar result as above, but in one direction) in \cite[Theorem 1.2.1]{Alpay} and \cite[Theorem 2.8]{Brod}. However, here out results are different in the following sense: (i) we are interested in scalar-valued (unlike operator-valued functions in \cite{Alpay, Brod}) Schur functions, (ii) our isometric colligations are explicit, (iii) our method is reversible (see Subsection \ref{sub-reverse}), and (perhaps most importantly) (iv) our ideas works in the setting of $n$-variable Schur(-Agler) functions.

We continue the discussion by presenting a transfer function realization of a two variables Schur function (see \cite{AM1} and also page 171, \cite{AM}):

\begin{thm}[Agler]\label{thm-Agler D2}
Let $\vp$ be a function on $\D^2$. Then $\vp \in \cls (\D^2)$ if and only if there exist Hilbert spaces $\clh_1$ and $\clh_2$ and an isometric colligation
\[
V = \begin{bmatrix} a & B \\ C & D \end{bmatrix} \in \clb(\mathbb{C} \oplus (\clh_1 \oplus \clh_2)),
\]
such that $\vp = \tau_V$ where
\[
\tau_V(\z) = a + B (I_{\clh_1 \oplus \clh_2} -E_{\clh_1 \oplus \clh_2}(\z) D)^{-1} E_{\clh_1 \oplus \clh_2}(\z) C,
\]
and $E_{\clh_1 \oplus \clh_2}(\z) = z_1 I_{\clh_1} \oplus z_2 I_{\clh_2}$ for all $\z \in \D^2$.
\end{thm}
Here and throughout the paper, elements of $\mathbb{C}^n$ will be denoted by $\z$, that is, $\z = (z_1, \ldots, z_n) \in \mathbb{C}^n$. Also we denote by $\clb(\clh_1, \clh_2)$ (and simply by $\clb(\clh_1)$ if $\clh_1 = \clh_2$) the set of all bounded linear operators from the Hilbert space $\clh_1$ into the Hilbert space $\clh_2$.

Agler's result exemplify the possibility of transfer function realizations (corresponding to isometric colligations) of Schur functions in $n$-variables, $n > 2$. This is, however, not true in general, and the possibility of transfer function realizations of functions in $\cls(\D^n)$, $n \geq 3$, is closely related to (as also the ideas in Agler's proof suggests) the subtlety of von Neumann inequality of commuting $n$-tuples of contractions, $n > 2$, on Hilbert spaces.

This motivates consideration of a special class of bounded analytic functions: The \textit{Schur-Agler class} $\cls \cla(\D^n)$ \cite{JA1} consists of scalar-valued analytic functions $\vp$ on $\D^n$ such that $\vp$ satisfies the $n$-variables von Neumann inequality, that is
\[
\|\vp(T_1, \ldots, T_n)\|_{\clb(\clh)} \leq 1,
\]
for any $n$-tuples of commuting strict contractions on a Hilbert space $\clh$. The elements of $\cls \cla(\D^n)$ are called \textit{Schur-Agler class functions}. If $\vp \in \cls \cla(\D^n)$, then we also say that $\vp$ is a \textit{function in the Schur-Agler class} $\cls \cla(\D^n)$. The following theorem of Jim Agler \cite{JA1} then obtains:

\begin{thm}[Agler]\label{thm-Agler Dn}
Let $\vp$ be a function on $\D^n$. Then $\vp \in \cls \cla(\D^n)$ if and only if there exist Hilbert spaces $\clh_1, \ldots, \clh_n$ and an isometric colligation
\[
V = \begin{bmatrix} a & B \\ C & D \end{bmatrix} \in \clb(\mathbb{C} \oplus \clh_1^n),
\]
such that $\vp = \tau_V$ where
\[
\tau_V(\z) = a + B  (I_{\clh_1^n} - E_{\clh_1^n}(\z) D)^{-1} E_{\clh_1^n}(\z) C,
\]
$\clh_1^n = \displaystyle \bigoplus_{i=1}^n \clh_i$ and $E_{\clh_1^n}(\z) = \displaystyle \bigoplus_{i=1}^n z_i I_{\clh_i}$ for all $\z \in \D^n$.
\end{thm}

Following the classical (one variable) von Neumann inequality, Ando \cite{Ando} proved that the von Neumann inequality also holds for commuting pairs of contractions. On the other hand, as we have pointed out earlier, the von Neumann inequality does not hold in general for $n$-tuples, $n > 2$, of commuting contractions \cite{CD, Var}. It follows then that
\[
\cls(\D) = \cls \cla(\D) \quad \mbox{and} \quad \cls(\D^2) = \cls \cla(\D^2),
\]
but $\cls(\D^n) \supsetneq \cls \cla(\D^n)$ for all $n > 2$.

Needless to say, transfer function realizations and isometric colligation matrices corresponding to Schur-Agler class functions in $n$-variables, $n \geq 1$, are among the most frequently used techniques in problems in function theory, operator theory and interdisciplinary subjects such as Nevanlinna-Pick interpolation \cite{AM1}, commutant lifting theorem and analytic model theory \cite{DS, F, FFGK}, scattering theory \cite{AG}, interpolation and Toeplitz corona theorem \cite{BTV}, electrical network theory \cite{Hel, Helton}, signal processing \cite{Kailath, IG}, linear systems \cite{Kaczorek, DFT, T}, operator algebras \cite{MP, MS} and image processing \cite{Roesser} (just to name a few). In this context and for deeper studies, we refer the reader to a number of classic work such as Liv\v{s}ic \cite{L1, L2}, Brodski\u{i} \cite{Brod}, Brodski\u{i} and M. Liv\v{s}ic \cite{BL1} and Pavlov \cite{Pav}. Also see \cite{ADR}, \cite{BGEB} and \cite{GVVW} and the references therein.

From this point of view, along with a question of interest in its own right, here we aim at finding necessary and sufficient conditions on isometric colligations which guarantee that a Schur-Agler class function factors into a product of Schur-Agler class functions. More precisely, we aim to solve the following problem: Given $\theta \in \cls \cla(\D^n)$, find a set of necessary and sufficient conditions on isometric colligations $V$ which ensures that
\[
\theta = \tau_V = \phi \psi,
\]
for some (explicit) $\phi$ and $\psi$ in $\cls \cla(\D^n)$.

In this paper we give a complete answer to this question by identifying checkable conditions on isometric colligations. Our results and approach are new even in the case of one variable and two-variable Schur functions (however, see the paragraph preceding Theorem \ref{thm-Agler D2}). In this context, it is also worth noting that the structure of bounded analytic functions in several variables is much more complicated than the structure of Schur functions on the unit disc (for instance, consider the existence of inner-outer factorizations of bounded analytic functions in one variable). From this point of view, our approach is also focused on providing an understanding of the complex area of bounded analytic functions of two or more variables (as the transfer function realization technique has already proven to be extremely useful in proving many classical results like Nevanlinna-Pick interpolation theorem and Carath\'{e}odory interpolation theorem etc. in several variables).

Our main results, specializing to the $n=2$ case, yields the following: Suppose $\theta \in \cls(\D^2)$ and $a :=\theta(0) \neq 0$. Then:

\NI (1) \textit{Theorem \ref{th-factor 2} implies that}: $\theta(\z) = \vp_1(z_1) \vp_2(z_2)$, $\z \in \D^2$, for some $\vp_1$ and $\vp_2$ in $\cls(\D)$ if and only if $\theta = \tau_V$ for some isometric colligation
\[
V = \left[\begin{array}{@{}c|cc@{}}
a & B_1 & B_2
\\ \hline
C_1 & D_{11} & \frac{1}{a} C_1 B_2 \\
C_2 & 0 & D_{22}
\end{array}\right] \in \clb(\mathbb{C} \oplus (\clh_1 \oplus \clh_2)).
\]

\NI \textit{(2) Theorem \ref{th-factor 3} implies that}: $\theta = \vp \psi$ for some $\vp$ and $\psi$ in $\cls(\D^2)$ if and only if there exist Hilbert spaces $\{\clm_i\}_{i=1}^2$ and $\{\cln_i\}_{i=1}^2$ and isometric colligation
\[
V =  \left[\begin{array}{@{}c|cc@{}}
a & B_1 & B_2
\\ \hline
C_1 & D_{11} & D_{12} \\
C_2 & D_{21} & D_{22}
\end{array}\right] \in \clb\Big(\mathbb{C} \oplus ((\clm_1 \oplus \cln_1) \oplus (\clm_2 \oplus \cln_2))\Big),
\]
such that $\theta = \tau_V$, and representing $B_i$, $C_i$ and $D_{ij}$ as
\[
B_i = \begin{bmatrix} B_i(1) & B_i(2) \end{bmatrix} \in \clb(\clm_i \oplus \cln_i, \mathbb{C}) \quad \mbox{and} \quad C_i = \begin{bmatrix} C_i(1) \\ C_i(2) \end{bmatrix} \in \clb(\mathbb{C}, \clm_i \oplus \cln_i),
\]
and $D_{ij} = \begin{bmatrix} D_{ij}(1) & D_{ij}(12) \\ D_{ij}(21) & D_{ij}(2) \end{bmatrix} \in \clb(\clm_j \oplus \cln_j, \clm_i \oplus \cln_i)$, respectively, one has $D_{ij}(21) = 0$ and $a D_{ij}(2) = C_i(1) B_j(2)$, $i, j = 1, 2$.

\NI Moreover, in the case of (1) (see Theorem \ref{thm-sufficient V1V2}): $\phi_1(z) = \tau_{\tilde{V}_1}(z)$ and $\phi_2(z) = \tau_{\tilde{V}_2}(z)$, $z \in \D$, where
\[
\tilde{V}_1 =  \begin{bmatrix}
\alpha & B_1 \\
\frac{1}{\beta} C_1 & D_{11} \end{bmatrix} \quad \mbox{and} \quad \tilde{V}_2 = \begin{bmatrix}
\beta & \frac{1}{\alpha}B_2
\\ C_2 & D_{22}
\end{bmatrix},
\]
and $\alpha$ and $\beta$ are non-zero scalars satisfying the conditions $|\beta|^2 = 1 - C_2^* C_2$ and $\alpha = \frac{a}{\beta}$; and in the case of (2) (see Theorem \ref{thm-3 between}): $\vp(\z) = \tau_{V_1}(\z)$ and $\psi(\z) = \tau_{V_2}(\z)$, $\z \in \D^2$, where
\[
V_1 = \begin{bmatrix} \alpha & B(1) \\ \frac{1}{\beta} C(1) & D(1) \end{bmatrix}
\quad \mbox{and} \quad
 V_2 = \begin{bmatrix} \beta & \frac{1}{\alpha} B(2) \\ C(2) & D(2) \end{bmatrix},
\]
and
\[ D(1) = \begin{bmatrix} D_{kl}(1)\end{bmatrix}_{k,l=1}^2,\;  D(2) = \begin{bmatrix} D_{kl}(2)\end{bmatrix}_{k,l=1}^2,\;
B(i) = \begin{bmatrix} B_1(i) & B_2(i) \end{bmatrix} \text{ and } C(i) = \begin{bmatrix} C_1(i) \\ C_2(i) \end{bmatrix},
\]
for all $i = 1, 2$,
and $\alpha$ and $\beta$ are non-zero scalars satisfying the conditions  $|\beta|^2 = |a|^2 + C(1)^* C(1)$ and $\alpha = \frac{a}{\beta}$.

\begin{rem}\label{rem-1}
The assumption that $\theta(0) \neq 0$ is not essential for the necessary parts of the above results (and Theorems \ref{th-factor 2} and \ref{th-factor 3}) and the case of $\theta(0) = 0$ will be treated separately in Section \ref{sbusect theta 0}. As we will see there, functions vanishing at the origin reveals more detailed properties of corresponding isometric colligations.
\end{rem}

The rest of this paper is organized as follows. Section \ref{sec-2} contains the definition of $\clf_m(n)$ class of isometric colligations, $1 \leq m < n$, and a classification of factorizations of functions in the Schur-Agler class $\cls\cla(\D^n)$, $n > 1$, into Schur-Agler class factors with fewer variables. Section \ref{sec-factor Fn} introduces the $\clf(n)$ class of isometric colligations, which connects the representation of a Schur-Agler class function to its Schur-Agler class factors.  In Section \ref{sbusect theta 0}, we will discuss factorizations of Schur-Agler class functions vanishing at the origin. The concluding Section \ref{sec last} outlines some concrete examples and presents results concerning one variable factors of Schur-Agler class functions and a remark on the reversibility of our method of factorizations.

\newsection{Factorizations and Property $\clf_m(n)$}\label{sec-2}

In this section, we present results concerning factorizations of Schur-Agler class functions in $\cls\cla(\D^n)$, $n > 1$, into Schur-Agler class factors with fewer variables. More specifically, our interest here is to identify (and then classify) isometric colligations $V$ such that $\tau_V \in \cls\cla (\D^n)$ and
\[
\tau_V(\z) = \phi(z_1, \ldots, z_m) \psi(z_{m+1}, \ldots, z_n) \quad\quad (\z \in \D^n),
\]
for some (canonical, in terms of $V$) $\phi \in \cls\cla(\D^m)$ and $\psi \in \cls\cla(\D^{n-m})$. Throughout this section we will always assume that $1 \leq m < n$.

We begin with fixing some notation. Given $1 \leq m < p \leq n$ and Hilbert spaces $\clh_1, \ldots, \clh_n$, we set
\[
\clh_m^p = \clh_m \oplus \clh_{m+1} \oplus \cdots \oplus \clh_p.
\]
In particular, $\clh_1^n = \bigoplus\limits_{i=1}^n \clh_i$. Moreover, with respect to the orthogonal decomposition $\clh_1^n = \clh_1^m \oplus \clh_{m+1}^n$, we represent an operator $D \in \clb(\clh_1^n)$ as
\[
D = \begin{bmatrix} D_{11} & D_{12}\\ D_{21} & D_{22} \end{bmatrix} \in \clb(\clh_1^m \oplus \clh_{m+1}^n).
\]
Similarly, if $\cle$ and $\cle_*$ are Hilbert spaces, $B \in \clb(\clh_1^n, \cle)$ and $C \in \clb(\cle_*, \clh_1^n)$, then we write
\[
B = \begin{bmatrix} B_1 & B_2 \end{bmatrix} \in \clb(\clh_1^m \oplus \clh_{m+1}^n, \cle) \quad \mbox{and} \quad
C = \begin{bmatrix} C_1 \\ C_2 \end{bmatrix} \in \clb(\cle_*, \clh_1^m \oplus \clh_{m+1}^n).
\]
Now we are ready to introduce the central object of this section.

\begin{defn}
Let $1 \leq m < n$. We say that an isometry $V \in \clb(\clh)$ satisfies property $\clf_m(n)$ if there exist Hilbert spaces $\clh_1, \ldots, \clh_n$ such that $\clh = \mathbb{C} \oplus \clh_1^n$, and representing $V$ as
\[
V = \begin{bmatrix} a & B_{1} & B_{2} \\ C_1 & D_{11} & D_{12}\\ C_2 & D_{21} & D_{22} \end{bmatrix} \in \clb(\mathbb{C} \oplus \clh_1^m \oplus \clh_{m+1}^n),
\]
one has $D_{21} = 0$ and $a D_{12} = C_1 B_2$.
\end{defn}

More specifically, an isometry $V \in \clb(\clh)$ satisfies property $\clf_m(n)$ if there exist Hilbert spaces $\clh_1, \ldots, \clh_n$ such that $\clh = \mathbb{C}\oplus \mathcal{H}_{1}\oplus\cdots \oplus\mathcal{H}_{n}$, and writing $V$ as
\[
V=\left[\begin{array}{@{}c|ccc@{}}
a & B_{1}  & \cdots & B_{n} \\\hline
C_{1} & D_{11} & \cdots & D_{1 n} \\
\vdots & \vdots  & \ddots & \vdots \\
C_{n} & D_{n 1} & \cdots  & D_{n n}
\end{array}\right],
\]
on $\mathbb{C}\oplus (\mathcal{H}_{1}\oplus\cdots \oplus\mathcal{H}_{n})$, one has
\[
D_{ij}=0,
\]
for all $i= m + 1, \cdots, n$ and $j=1,\cdots, m$, and
\[
a D_{ij}=C_{i}B_{j},
\]
for all $i=1,\cdots, m$ and $j = m+1, \cdots, n$. By way of example, we consider the two variables situation. We say that an isometry $V$ satisfies property $\clf_1(2)$ if there exist Hilbert spaces $\clh_1$ and $\clh_2$ such that
\[
V =  \begin{bmatrix} a & B_{1} & B_{2} \\ C_1 & D_{11} & D_{12}\\ C_2 & 0 & D_{22} \end{bmatrix} \in \clb(\mathbb{C} \oplus \clh_1 \oplus \clh_2)
\]
and $a D_{12} = C_1 B_2$.

Let us introduce some more notation. Let $1 \leq m < p \leq n$. We set
\[
E_{\clh_m^p}(\z) = z_m I_{\clh_m} \oplus \cdots \oplus z_p I_{\clh_p} \quad \quad (\z \in \mathbb{C}^n).
\]
Also for $X \in \clb(\clh_m^p)$, $\|X\| \leq 1$, define
\[
R_m^p(\z, X) = \Big(I_{\clh_m^p} - E_{\clh_m^p}(\z) X \Big)^{-1} \quad \quad (\z \in \D^n).
\]
Note that $R_m^p(\z, X)$ is a function of $\{z_m, \ldots, z_p\}$ variables. Moreover, we will denote $R_1^n(\z, X)$ simply by $R(\z, X)$.

Now we proceed to prove that a pair of isometric colligations is naturally associated with an isometric colligation satisfying property $\clf_m(n)$. More specifically, given $\tau_{V_1} \in \cls \cla(\D^m)$ and $\tau_{V_2} \in \cls \cla (\D^{n-m})$ for some isometric colligations $V_1$ and $V_2$, we aim to construct an explicit isometric colligation $V$ such that $V$ satisfies property $\clf_m(n)$ and
\[
\tau_V(\z) = \tau_{V_1}(z_1, \ldots, z_m) \tau_{V_2}(z_{m+1}, \ldots, z_n) \quad \quad (\z \in \D^n).
\]
To this end, let $\clh_1, \ldots, \clh_n$ be Hilbert spaces. Suppose
\[
V_1 = \begin{bmatrix}a_1 & B_1 \\ C_1 & D_1 \end{bmatrix} \in \clb(\mathbb{C} \oplus \clh_1^m),
\quad \mbox{and}\quad
V_2 = \begin{bmatrix}a_2 & B_2 \\ C_2 & D_2 \end{bmatrix} \in \clb(\mathbb{C} \oplus \clh_{m+1}^n),
\]
are isometric colligations. Define $\tilde{V}_1$ and $\tilde{V}_2$ in $\clb(\mathbb{C} \oplus \clh_1^m \oplus \clh_{m+1}^n)$ by
\[
\tilde{V}_1 = \begin{bmatrix}a_1 & B_1 & 0 \\ C_1 & D_1 & 0 \\ 0 & 0 & I \end{bmatrix}
\quad \mbox{and} \quad
\tilde{V}_2 = \begin{bmatrix}a_2 & 0 & B_2 \\ 0 & I & 0 \\ C_2 & 0 & D_2 \end{bmatrix},
\]
and set $V = \tilde{V}_1 \tilde{V}_2$. It is easy to check, by swapping rows and columns (of $\tilde{V}_2$), that $\tilde{V}_1$ and $\tilde{V}_2$ are isometries and thus the isometric colligation
\[
V = \left[\begin{array}{@{}c|cc@{}}
a_1 a_2 & B_1 & a_1 B_2
\\ \hline
a_2 C_1 & D_1 & C_1 B_2 \\ C_2 & 0 & D_2
\end{array}\right] \in \clb\Big(\mathbb{C} \oplus (\clh_1^m \oplus \clh_{m+1}^n)\Big),
\]
satisfies property $\clf_m(n)$. Let $\z \in \D^n$. Clearly
\[
\tau_V(\z) = a_1 a_2 + \begin{bmatrix} B_1 & a_1 B_2 \end{bmatrix} R(\z, \begin{bmatrix} D_1 & C_1 B_2 \\ 0 & D_2 \end{bmatrix}) E_{\clh_1^n}(\z) \begin{bmatrix} a_2 C_1 \\ C_2 \end{bmatrix},
\]
where
\[
\begin{split}
R\Big(\z, \begin{bmatrix} D_1 & C_1 B_2 \\ 0 & D_2 \end{bmatrix}\Big)^{-1} & = I_{\clh_1^n} - \begin{bmatrix}E_{\clh_1^m}(\z) & 0 \\ 0 & E_{\clh_{m+1}^n}(\z) \end{bmatrix} \begin{bmatrix} D_1 & C_1 B_2 \\ 0 & D_2 \end{bmatrix}
\\
& = \begin{bmatrix}I_{\clh_1^m} - E_{\clh_1^m}(\z) D_1 & - E_{\clh_1^m}(\z) C_1 B_2 \\ 0 & I_{\clh_{m+1}^n} - E_{\clh_{m+1}^n}(\z) D_2 \end{bmatrix}.
\end{split}
\]
By the inverse formula of an invertible upper triangular matrix, it follows that
\[
R\Big(\z, \begin{bmatrix} D_1 & C_1 B_2 \\ 0 & D_2 \end{bmatrix}\Big) = \begin{bmatrix} R_1^m(\z, D_1) & R_1^m(\z, D_1) E_{\clh_1^m}(\z) C_1 B_2 R_{m+1}^n(\z, D_2) \\ 0 & R_{m+1}^n(\z, D_2) \end{bmatrix}.
\]
We now infer, in view of the above equality, that
\[
\begin{split}
\tau_V(\z) & = a_1 a_2 + \begin{bmatrix} B_1 & a_1 B_2 \end{bmatrix} R\Big(\z, \begin{bmatrix} D_1 & C_1 B_2 \\ 0 & D_2 \end{bmatrix}\Big) E_{\clh_1^n}(\z) \begin{bmatrix} a_2 C_1 \\ C_2 \end{bmatrix}
\\
& = a_1 a_2 + \begin{bmatrix} B_1 & a_1 B_2 \end{bmatrix} \begin{bmatrix} R_1^m(\z, D_1) & R_1^m(\z, D_1) E_{\clh_1^m}(\z) C_1 B_2 R_{m+1}^n(\z, D_2) \\ 0 & R_{m+1}^n(\z, D_2) \end{bmatrix}
\\
& \quad  \quad \quad \times \begin{bmatrix} a_2 E_{\clh_1^m}(\z) C_1 \\ E_{\clh_{m+1}^n}(\z) C_2 \end{bmatrix}
\\
& = a_1 a_2 + a_2 B_1 R_1^m(\z, D_1)  E_{\clh_1^m}(\z) C_1  + a_1 B_2 R_{m+1}^n(\z, D_2) E_{\clh_{m+1}^n}(\z) C_2
\\
& \quad + B_1 R_1^m(\z, D_1) E_{\clh_1^m}(\z) C_1 B_2 R_{m+1}^n(\z, D_2)
E_{\clh_{m+1}^n}(\z) C_2
\\
& = \Big(a_1 + B_1 R_1^m(\z, D_1)  E_{\clh_1^m}(\z) C_1\Big) \Big(a_2 + B_2 R_{m+1}^n (\z, D_2)  E_{\clh_{m+1}^n}(\z) C_2 \Big)
\\
& = \tau_{V_1}(z_1, \ldots, z_m) \tau_{V_2}(z_{m+1}, \ldots, z_n),
\end{split}
\]
for all $\z \in \D^n$. We have therefore proved the following result:

\begin{thm}\label{th-1st factor}
Let $1 \leq m <n$, and let $\clh_1, \ldots, \clh_n$ be Hilbert spaces. Suppose
\[
V_1 = \begin{bmatrix}a_1 & B_1 \\ C_1 & D_1 \end{bmatrix} \in \clb(\mathbb{C} \oplus (\displaystyle \bigoplus_{i=1}^m \clh_i)) \quad \mbox{and} \quad V_2 = \begin{bmatrix}a_2 & B_2 \\ C_2 & D_2 \end{bmatrix} \in \clb(\mathbb{C} \oplus  (\displaystyle \bigoplus_{i=m+1}^n \clh_i)),
\]
are isometric colligations. Define $\tilde{V}_1$, $\tilde{V}_2$ and $V$ in $\clb\Big(\mathbb{C} \oplus \Big((\displaystyle \bigoplus_{i=1}^m \clh_i) \oplus (\displaystyle \bigoplus_{i=m+1}^n \clh_i)\Big)\Big)$ by
\[
\tilde{V}_1 = \begin{bmatrix}a_1 & B_1 & 0 \\ C_1 & D_1 & 0 \\ 0 & 0 & I \end{bmatrix} \quad \mbox{and} \quad \tilde{V}_2 = \begin{bmatrix}a_2 & 0 & B_2 \\ 0 & I & 0 \\ C_2 & 0 & D_2 \end{bmatrix},
\]
and $V = \tilde{V}_1 \tilde{V}_2$, respectively. Then
\[
V =\left[\begin{array}{@{}c|cc@{}}
a_1 a_2 & B_1 & a_1 B_2 \\\hline
a_2 C_1 & D_1 & C_1 B_2 \\
C_2 & 0 & D_2
\end{array}\right] \in \clb\Big(\mathbb{C} \oplus \Big((\displaystyle \bigoplus_{i=1}^m \clh_i) \oplus (\displaystyle \bigoplus_{i=m+1}^n \clh_i)\Big)\Big),
\]
is an isometric colligation, $V$ satisfies property $\clf_m(n)$ and
\[
\tau_V(\z) = \tau_{V_1}(z_1, \ldots, z_m) \tau_{V_2}(z_{m+1}, \ldots, z_n) \quad \quad (\z \in \D^n).
\]
\end{thm}

Now to prove the reverse direction, \textit{we assume in addition that $\tau_V(0) \neq 0$} (for the case of transfer functions vanishing at the origin, see Section \ref{sbusect theta 0}) : Suppose $\clh_1, \ldots, \clh_n$ are Hilbert spaces and
\begin{equation}\label{eq-V33}
V = \begin{bmatrix} a & B_1 & B_2 \\ C_1 & D_{11} & D_{12} \\ C_2 & 0 & D_{22} \end{bmatrix} \in \clb(\mathbb{C} \oplus \clh_1^m \oplus \clh_{m+1}^n),
\end{equation}
is an isometric colligation satisfying property $\clf_m(n)$. Thus
\begin{equation}\label{eq-aDCB}
a D_{12} = C_1 B_2.
\end{equation}
Suppose $a := \tau_V(0) \neq 0$. Since $V^* V = I$, we have
\[
|a|^2 + C_1^* C_1 + C_2^* C_2 = 1,
\]
implies that
\[
1 - C_2^* C_2 = |a|^2 + C_1^* C_1 >0,
\]
as $a \neq 0$. Then there exists a scalar $\beta$, $0 < |\beta| \leq 1$, such that
\[
|\beta|^2 = 1 - C_2^* C_2.
\]
It now follows that
\begin{equation}\label{eq-a2aC}
C_1^* C_1 = |\beta|^2 - |a|^2,
\end{equation}
and
\begin{equation}\label{eq-aaa}
\alpha := \frac{a}{\beta},
\end{equation}
is a non-zero scalar. Define
\[
V_1 = \begin{bmatrix} \alpha & B_1 & 0 \\ \frac{1}{\beta} C_1 & D_{11} & 0 \\ 0 & 0 & I \end{bmatrix} \quad \mbox{and} \quad V_2 = \begin{bmatrix} \beta & 0 & \frac{1}{\alpha}B_2 \\ 0& I & 0 \\ C_2 & 0 & D_{22} \end{bmatrix},
\]
on $\mathbb{C} \oplus \clh_1^m \oplus \clh_{m+1}^n$.
It follows from \eqref{eq-a2aC} and \eqref{eq-aaa} that
\[
|\alpha|^2 + \frac{1}{|\beta|^2} C_1^* C_1 = |\alpha|^2 + \frac{1}{|\beta|^2} (|\beta|^2 - |a|^2)
= 1 + |\alpha|^2 - \frac{|a|^2}{|\beta|^2}
 = 1
\]
that is
\begin{equation}\label{eq-a1a21}
|\alpha|^2 + \frac{1}{|\beta|^2} C_1^* C_1 = 1.
\end{equation}
Also,  we see that $B_1^* B_1 + D_{11}^* D_{11} = I$, and
\[
\bar{\alpha} B_1 + \frac{1}{\bar{\beta}} C_1^* D_{11}  = \frac{1}{\bar{\beta}}(\bar{\alpha} \bar{\beta} B_1 + C_1^* D_{11})
= \frac{1}{\bar{\beta}}(\bar{a} B_1 + C_1^* D_{11})= 0,
\]
and hence $V_1^* V_1 = I$. We now proceed to prove that $V_2$ is also an isometry. First, it easy to see that $\bar{a} B_2 + C_1^* D_{12} + C_2^* D_{22} = 0$, and hence, by \eqref{eq-aDCB}, we have
\[
0 = \bar{a} B_2 + C_1^* D_{12} + C_2^* D_{22}
= \bar{a} B_2 + \frac{1}{a} C_1^* C_1 B_2 + C_2^* D_{22}
 = \frac{\bar{a}_2}{\alpha} (|\alpha|^2 + \frac{1}{|\beta|^2} C_1^* C_1) B_2 + C_2^* D_{22}.
\]
Then \eqref{eq-a1a21} implies that $\frac{\bar{\beta}}{\alpha} B_2 + C_2^* D_{22} = 0.$
Finally, again from $V^*V=I$ we get
\[
B_2^* B_2 + D_{12}^* D_{12} + D_{22}^* D_{22} = I.
\]
Now again by \eqref{eq-aDCB} we have
\[
\begin{split}
B_2^* B_2 + D_{12}^* D_{12} + D_{22}^* D_{22} & = B_2^*(1 + \frac{1}{|a|^2} C_1^* C_1) B_2 + D_{22}^* D_{22}
\\
& = \frac{1}{|\alpha|^2} B_2^*(|\alpha|^2 + \frac{1}{|\beta|^2} C_1^* C_1) B_2 + D_{22}^* D_{22},
\end{split}
\]
so that $\frac{1}{|\alpha|^2} B_2^* B_2 + D_{22}^* D_{22} = I,$
by \eqref{eq-a1a21}, from which we conclude that $V_2^* V_2 = I$. Finally, notice that
\[
V_1 V_2 = \begin{bmatrix} \alpha \beta & B_1 & B_2
\\
C_1 & D_{11} & \frac{1}{\alpha \beta} C_1 B_2
\\
C_2 & 0 & D_{22}
\end{bmatrix} = \begin{bmatrix} a & B_1 & B_2
\\
C_1 & D_{11} & \frac{1}{a} C_1 B_2
\\
C_2 & 0 & D_{22}
\end{bmatrix},
\]
and hence $V = V_1 V_2$, by \eqref{eq-aDCB}. Then, by Theorem \ref{th-1st factor}, we have
\[
\tau_V(\z) = \tau_{\tilde{V}_1}(z_1, \ldots, z_m) \tau_{\tilde{V}_2}(z_{m+1}, \ldots, z_n),
\]
for all $\z \in \D^n$ where $\tilde{V}_1 = \begin{bmatrix} \alpha & B_1 \\ \frac{1}{\beta} C_1 & D_{11} \end{bmatrix}$ and $\tilde{V}_2 = \begin{bmatrix} \beta & \frac{1}{\alpha} B_2 \\  C_2 & D_{22} \end{bmatrix}$. Thus we have proved the following statement:

\begin{thm}\label{thm-sufficient V1V2}
Suppose $\clh_1, \ldots, \clh_n$ are Hilbert spaces and $a$ be a non-zero scalar. If
\[
V = \left[\begin{array}{@{}c|cc@{}}
a & B_1 & B_2 \\
\hline
C_1 & D_{11} & \frac{1}{a} C_1 B_2 \\ C_2 & 0 & D_{22}
\end{array}\right] \in \clb \Big(\mathbb{C} \oplus \Big((\displaystyle \bigoplus_{i=1}^m \clh_i ) \oplus (\displaystyle \bigoplus_{i=m+1}^n \clh_i)\Big)\Big),
\]
is an isometric colligation, then
\[
\tilde{V}_1 = \begin{bmatrix} \alpha & B_1 \\ \frac{1}{\beta} C_1 & D_{11} \end{bmatrix} \quad \mbox{and} \quad \tilde{V}_2 = \begin{bmatrix} \beta & \frac{1}{\alpha} B_2 \\  C_2 & D_{22} \end{bmatrix}.
\]
are isometric colligations in $\clb \Big(\mathbb{C} \oplus (\displaystyle \bigoplus_{i=1}^m \clh_i )\Big)$ and $\clb \Big(\mathbb{C} \oplus (\displaystyle \bigoplus_{i=m+1}^n \clh_i)\Big)$, respectively, and
\[
\tau_V(\z) = \tau_{\tilde{V}_1}(z_1, \ldots, z_m) \tau_{\tilde{V}_2}(z_{m+1}, \ldots, z_n) \quad \quad (\z \in \D^n),
\]
where $\alpha$ and $\beta$ are non-zero scalars and satisfies the following conditions
\[
|\beta|^2 = |a|^2 + C_1^* C_1 \quad \mbox{and} \quad \alpha = \frac{a}{\beta}.
\]
\end{thm}

Summing up the results of Theorems \ref{th-1st factor} and \ref{thm-sufficient V1V2}, we conclude the following factorization theorem on Schur-Agler class functions in $\cls \cla(\D^n)$, $n \geq 2$:

\begin{thm}\label{th-factor 2}
Let $1 \leq m < n$, and let $\theta\in \mathcal{SA}(\mathbb{D}^n)$. If $\theta(0)\neq 0$, then
\[
\theta(\z) = \phi(z_1,\cdots,z_m) \psi(z_{m+1},\cdots,z_n) \quad \quad (\z \in \D^n),
\]
for some $\phi \in \mathcal{SA}(\mathbb{D}^m)$ and $\psi\in \mathcal{SA}(\mathbb{D}^{n-m})$	if and only if
\[
\theta(\z) = \tau_V(\z) \quad \quad(\z \in \D^n),
\]
for some isometric colligation $V$ satisfying property $\clf_m(n)$.
\end{thm}

We again point out that the assumption $\theta(0) \neq 0$ is not needed to prove the necessary part of the above theorem. Classification of factorizations of functions vanishing at the origin will be discussed in detail in Section \ref{sbusect theta 0}.

\newsection{Factorizations and Property $\clf(n)$}\label{sec-factor Fn}

In this section we investigate general $n$-variables Schur-Agler class factors of  Schur-Agler class functions in $\cls \cla(D^n)$. More specifically, for a given $\theta \in \cls\cla(\D^n)$, we give a set of necessary and sufficient conditions on isometric colligations ensuring the existence of $\vp$ and $\psi$ in $\cls\cla(\D^n)$ such that  $\theta = \vp \psi$. We identify a new class of isometric colligations, namely $\clf(n)$, and prove that the (Schur-Agler class) factors of Schur-Agler class functions are completely determined by isometric colligations satisfying property $\clf(n)$. Here we do not set any restriction on $n$, that is, we will assume that $n \geq 1$.

We first identify the relevant isometric colligations:

\begin{defn}
We say that an isometry $V \in \clb(\clh)$ satisfies property $\clf(n)$ if there exist Hilbert spaces $\{\clm_i\}_{i=1}^n$ and $\{\cln_i\}_{i=1}^n$ such that
\[
\clh = \mathbb{C} \oplus \Big(\bigoplus_{i=1}^n (\clm_i \oplus \cln_i)\Big),
\]
and representing $V$ as
\[
V=\left[\begin{array}{@{}c|ccc@{}}
a &  B_{1}  & \cdots & B_{n} \\\hline
C_{1} & D_{11} & \cdots & D_{1n} \\
\vdots & \vdots  & \ddots & \vdots \\
C_{n} & D_{n1} & \cdots  & D_{nn}
\end{array}\right] \in \clb\Big(\mathbb{C} \oplus  \Big(\bigoplus_{i=1}^n (\clm_i \oplus \cln_i)\Big)\Big) ,
\]
and $B_i$, $C_i$ and $D_{ij}$ as
\[
B_i = \begin{bmatrix} B_i (1) & B_i (2) \end{bmatrix} \in \clb(\clm_i \oplus \cln_i, \mathbb{C}), \quad C_i = \begin{bmatrix} C_i (1) \\ C_i (2) \end{bmatrix} \in \clb(\mathbb{C}, \clm_i \oplus \cln_i),
\]
and
\[
D_{ij} = \begin{bmatrix}
D_{ij}(1)& D_{ij} (12)\\
D_{ij} (21) & D_{ij} (2)
\end{bmatrix} \in \clb(\clm_j \oplus \cln_j, \clm_i \oplus \cln_i),
\]
one has
\[
D_{ij}(21) = 0, \quad \mbox{and}\quad  a D_{ij} (12) = C_i(1) B_j(2),
\]
for all $i, j = 1, \ldots, n$.
\end{defn}

As in Section \ref{sec-2}, here we also first prove that a pair of isometric colligations is naturally associated with an isometric colligation satisfying property $\clf(n)$. Let $\{\clm_i\}_{i=1}^n$ and $\{\cln_i\}_{i=1}^n$ be Hilbert spaces, and let
\[
V_1 = \begin{bmatrix} \alpha & B \\ C & D \end{bmatrix} = \left[\begin{array}{@{}c|ccc@{}}
\alpha & B_{1}  & \cdots & B_{n} \\\hline
C_{1} & D_{11} & \cdots & D_{1n} \\
\vdots & \vdots  & \ddots & \vdots \\
C_{n} & D_{n1} & \cdots  & D_{nn}
\end{array}\right] \in \clb(\mathbb{C} \oplus \clm_1^n),
\]
and
\[
V_2 = \begin{bmatrix} \beta & F \\ G & H \end{bmatrix} = \left[\begin{array}{@{}c|ccc@{}}
\beta & F_{1}  & \cdots & F_{n} \\\hline
G_{1} & H_{11} & \cdots & H_{1n} \\
\vdots & \vdots  & \ddots & \vdots \\
G_{n} & H_{n1} & \cdots  & H_{nn}
\end{array}\right] \in \clb(\mathbb{C} \oplus \cln_1^n),
\]
be isometric colligations. Given $i =1, \ldots, n$, we define $\clh_i = \clm_i \oplus \cln_i$, and bounded linear operators $\tilde{B}_i$, $\tilde{C}_i$ and $\tilde{D}_{ij}$ as
\[
\tilde{B}_i = \begin{bmatrix} B_i & 0 \end{bmatrix} \in \clb(\clh_i, \mathbb{C}), \; \tilde{C}_i = \begin{bmatrix} C_i \\ 0 \end{bmatrix} \in \clb(\mathbb{C}, \clh_i), \;
\mbox{and}\; \tilde{D}_{ij} =
\begin{bmatrix} D_{ij} & 0 \\ 0 & \delta_{ij} I \end{bmatrix} \in \clb(\clh_j, \clh_i),
\]
for all $i, j = 1, \ldots, n$. Set
\begin{equation}\label{eq-tildeV1}
\tilde{V}_1= \left[\begin{array}{@{}c|ccc@{}}
\alpha & \tilde{B}_{1}  & \cdots & \tilde{B}_{n} \\\hline
\tilde{C}_{1} & \tilde{D}_{11} & \cdots & \tilde{D}_{1n} \\
\vdots & \vdots  & \ddots & \vdots \\
\tilde{C}_{n} & \tilde{D}_{n1} & \cdots  & \tilde{D}_{nn}
\end{array}\right] \in \clb(\mathbb{C} \oplus \clh_1^n).
\end{equation}
On the other hand, let
\begin{equation}\label{eq-tildeV2}
\tilde{V}_2 = \left[\begin{array}{@{}c|ccc@{}}
\beta & \tilde{F}_{1}  & \cdots & \tilde{F}_{n} \\ \hline
\tilde{G}_{1} & \tilde{H}_{11} & \cdots & \tilde{H}_{1n} \\
\vdots & \vdots  & \ddots & \vdots \\
\tilde{G}_{n} & \tilde{H}_{n1} & \cdots & \tilde{H}_{nn}
\end{array}\right] \in \clb(\mathbb{C} \oplus \clh_1^n),
\end{equation}
where
\[
\tilde{F}_i = \begin{bmatrix} 0 & F_i \end{bmatrix} \in \clb(\clh_i, \mathbb{C}), \; \tilde{G}_i =  \begin{bmatrix} 0 \\ G_i \end{bmatrix} \in \clb(\mathbb{C}, \clh_i),
\;\mbox{and}\;
\tilde{H}_{ij} =
\begin{bmatrix} \delta_{ij} I & 0 \\ 0 & H_{ij} \end{bmatrix} \in \clb(\clh_j, \clh_i),
\]
for all $i, j = 1, \ldots, n$. Define $V = \tilde{V}_1 \tilde{V_2}$. It then follows that $V \in \clb(\mathbb{C} \oplus \clh_1^n)$ is an isometry and
\begin{equation}\label{eq-Def V}
V = \left[\begin{array}{@{}c|ccc@{}}
\alpha \beta & \hat{B}_{1}  & \cdots & \hat{B}_{n} \\\hline
\hat{C}_{1} & \hat{D}_{11} & \cdots & \hat{D}_{1n} \\
\vdots & \vdots  & \ddots & \vdots \\
\hat{C}_{n} & \hat{D}_{n1} & \cdots  & \hat{D}_{nn}
\end{array}\right] := \begin{bmatrix} \alpha \beta & \hat{B} \\ \hat{C} & \hat{D} \end{bmatrix} ,
\end{equation}
where
\begin{equation}\label{eq-3.1 one}
\hat{B}_i = \begin{bmatrix} B_i & \alpha F_i \end{bmatrix} \in \clb(\clh_i, \mathbb{C}), \; \hat{C}_i = \begin{bmatrix} \beta C_i \\ G_i \end{bmatrix} \in \clb(\mathbb{C}, \clh_i),
\text{ and }
\hat{D}_{ij} =
\begin{bmatrix} D_{ij} & C_i F_j \\ 0 & H_{ij} \end{bmatrix},
\end{equation}
for all $i, j = 1, \ldots, n$. Define $X(\z) : \mathbb{C} \raro \mathbb{C}$, $\z \in \D^n$, by
\[
X(\z) = \hat{B} (I_{\clh_1^n} - E_{\clh}(\z) \hat{D})^{-1} E_{\clh}(\z) \hat{C}.
\]
Then $\tau_V(\z) = \alpha \beta + X(\z)$, $\z \in \D^n$. Next, define the flip operator $\eta : \clh_1^n \raro \clm_1^n \oplus \cln_1^n,$
by
\begin{equation}\label{eq-flip}
\eta \Big(\bigoplus_{i=1}^n (f_i \oplus g_i)\Big) = (\bigoplus_{i=1}^n f_i) \oplus (\bigoplus_{i=1}^n g_i),
\end{equation}
for all $f_i \in \clm_i$ and $g_i \in \cln_i$, $i = 1, \ldots, n$. Then $\eta$ is a unitary operator and so
\[
X(\z) = (\hat{B} \eta^*) \Big(I_{\clm_1^n \oplus \cln_1^n} - (\eta E_{\clh_1^n}(\z) \eta^*) (\eta \hat{D} \eta^*)\Big)^{-1} (\eta E_{\clh_1^n}(\z) \eta^*) (\eta \hat{C}).
\]
On the other hand, the definition of the flip operator $\eta$ reveals that
\[
\hat{B} \eta^* = \begin{bmatrix} B & \alpha F \end{bmatrix}, \; \eta \hat{C} = \begin{bmatrix} \beta C \\ G \end{bmatrix}, \;
\eta \hat{D} \eta^* = \begin{bmatrix} D & C F \\ 0 & H \end{bmatrix}, \;
\text{and} \;
\eta E_{\clh_1^n}(\z) \eta^* = \begin{bmatrix}  E_{\clm_1^n}(\z) & 0 \\ 0 &  E_{\cln_1^n}(\z)\end{bmatrix}.
\]
In particular, this yields
\[
I_{\clm_1^n \oplus \cln_1^n} - (\eta E_{\clh_1^n}(\z) \eta^*) (\eta \hat{D} \eta^*) = \begin{bmatrix} I - E_{\clm_1^n}(\z)D & - E_{\clm_1^n}(\z) CF \\ 0 & I - E_{\cln_1^n}(\z) H \end{bmatrix} \qquad (\z \in \D^n).
\]
In order to further ease the notation, for Hilbert spaces $\{\cls_i\}_{i=1}^n$ and $\z \in \D^n$, we set
\[
E_{\cls}(\z) = \bigoplus_{i=1}^n z_i I_{\cls_i},
\]
and, for $Y \in \clb(\displaystyle \bigoplus_{i=1}^n \cls_i)$, $\|Y\| \leq 1$, define
$r(\z, Y) = \Big(I_{\cls_1^n} - E_{\cls}(\z) Y \Big)^{-1}.$

\NI Continuing the above computation, for each $\z \in \D^n$, we now have
\[
\Big(I_{\clm_1^n \oplus \cln_1^n} - (\eta E_{\clh_1^n}(\z) \eta^*) (\eta \hat{D} \eta^*)\Big)^{-1} = \begin{bmatrix} r(\z, D) & r(\z, D) E_{\clm_1^n}(\z) CF r(\z, H) \\ 0 & r(\z, H) \end{bmatrix}.
\]
Moreover, since $(\eta E_{\clh_1^n}(\z) \eta^*) (\eta \hat{C}) = \begin{bmatrix} \beta E_{\clm_1^n}(\z) C \\ E_{\cln_1^n}(\z) G \end{bmatrix}$, it follows that
\[
\begin{split}
X(\z) & = \begin{bmatrix} B & \alpha F \end{bmatrix} \begin{bmatrix} \beta r(\z, D) E_{\clm_1^n}(\z) C + r(\z, D) E_{\clm_1^n}(\z) C F r(\z, H) E_{\cln_1^n}(\z) G \\ r(\z, H)  E_{\cln_1^n}(\z) G  \end{bmatrix}
\\
& = \beta B r(\z, D) E_{\clm_1^n}(\z) C + B r(\z, D) E_{\clm_1^n}(\z) C F r(\z, H) E_{\cln_1^n}(\z) G + \alpha F r(\z, H)  E_{\cln_1^n}(\z) G,
\end{split}
\]
and so $\tau_V(\z) = \tau_{V_1}(\z) \tau_{V_2}(\z)$, $\z \in \D^n$. We have therefore proved:

\begin{thm}\label{th-MNH factor}
Suppose $V_1 = \begin{bmatrix} \alpha & B \\ C & D \end{bmatrix} \in \clb \Big(\mathbb{C} \oplus (\displaystyle \bigoplus_{i=1}^n \clm_i)\Big)$ and $V_2 = \begin{bmatrix} \beta & F \\ G & H \end{bmatrix} \in \clb\Big(\mathbb{C} \oplus  (\displaystyle \bigoplus_{i=1}^n \cln_i)\Big)$ are isometric colligations, and let $V = \tilde{V}_1 \tilde{V}_2$, where $\tilde{V}_1$ and $\tilde{V}_2$ are as in \eqref{eq-tildeV1} and \eqref{eq-tildeV2}, respectively.
Then the isometric colligation $V \in \clb\Big(\mathbb{C} \oplus (\displaystyle \bigoplus_{i=1}^n (\clm_i \oplus \cln_i))\Big)$ as in \eqref{eq-Def V} satisfies property $\clf(n)$ and $\tau_V = \tau_{V_1} \tau_{V_2}$.
\end{thm}

We have the following interpretations of the above theorem: Let $\theta, \phi, \psi \in \mathcal{SA}(\mathbb{D}^n)$, and suppose $\theta = \phi \psi$. Suppose $V_1 = \begin{bmatrix} \alpha & B \\ C & D \end{bmatrix}$ and $V_2 = \begin{bmatrix} \beta & F \\ G & H \end{bmatrix}$
are isometric colligations on $\mathbb{C} \oplus \clm_{1}^n$ and $\mathbb{C} \oplus  \cln_{1}^n$, respectively, and $\phi = \tau_{V_1}$, and $\psi = \tau_{V_2}$. Then the isometric colligation $V = \tilde{V}_1 \tilde{V}_2$, as constructed in Theorem \ref{th-MNH factor}, satisfies property $\clf(n)$ and $\tau_V(\z) = \tau_{V_1}(\z) \tau_{V_2}(\z)$ for all $\z \in \D^n$, that is, $\theta = \tau_V$.

Now we proceed to treat the converse of Theorem \ref{th-MNH factor}. Let $V \in \clb(\clh)$ be an isometric colligation, and let $V$ satisfies property $\clf(n)$. As in Theorem \ref{thm-sufficient V1V2}, here also we assume that $a := \tau_V(0) \neq 0$. Now
\[
\clh = \mathbb{C} \oplus \Big(\displaystyle \bigoplus_{i=1}^n (\clm_i \oplus\cln_i)\Big),
\]
for some Hilbert spaces $\{\clm_i\}_{i=1}^n$ and $\{\cln_i\}_{i=1}^n$, and
\begin{equation}\label{eq3 - Vblock}
V= \begin{bmatrix}a & B \\ C & D\end{bmatrix} = \left[\begin{array}{@{}c|ccc@{}}
a &  B_{1}  & \cdots & B_{n} \\\hline
C_{1} & D_{11} & \cdots & D_{1n} \\
\vdots & \vdots  & \ddots & \vdots \\
C_{n} & D_{n1} & \cdots  & D_{nn}
\end{array}\right],
\end{equation}
where
\begin{equation}\label{eq3 - BC}
B_i = \begin{bmatrix} B_i (1) & B_i (2) \end{bmatrix} \in \clb(\clm_i \oplus \cln_i, \mathbb{C}), \quad C_i = \begin{bmatrix} C_i (1) \\ C_i (2) \end{bmatrix} \in \clb(\mathbb{C}, \clm_i \oplus \cln_i),
\end{equation}
and
\begin{equation}\label{eq3 - D}
D_{ij} = \begin{bmatrix}
D_{ij}(1)& \frac{1}{a} C_i(1) B_j(2)\\
0 & D_{ij} (2)
\end{bmatrix} \in \clb(\clm_j \oplus \cln_j, \clm_i \oplus \cln_i),
\end{equation}
for all $i, j = 1, \ldots, n$. Set
\begin{equation}\label{eq3 - D1D2}
D(1) = \begin{bmatrix} D_{ij}(1) \end{bmatrix}_{i,j=1}^n \in \clb(\bigoplus_{i=1}^n \clm_i), \quad D(2) = \begin{bmatrix} D_{ij}(2) \end{bmatrix}_{i,j=1}^n \in \clb(\bigoplus_{i=1}^n \cln_i)
\end{equation}
and
\[
D(12) = \begin{bmatrix} D_{ij}(12) \end{bmatrix}_{i,j=1}^n \in \clb\Big(\bigoplus_{i=1}^n \clm_i, \bigoplus_{i=1}^n \cln_i \Big),
\]
and consider the flip operator $\eta : \Big(\bigoplus_{i=1}^n (\clm_i \oplus \cln_i)\Big) \raro (\bigoplus_{i=1}^n \clm_i) \oplus (\bigoplus_{i=1}^n \cln_i)$ (see \eqref{eq-flip}). Then
\[
\eta D \eta^* = \begin{bmatrix} D(1) & D(12) \\ 0 & D(2) \end{bmatrix} \in \clb\Big((\bigoplus_{i=1}^n \clm_i) \oplus (\bigoplus_{i=1}^n \cln_i)\Big).
\]
If we define $V_{\eta} := \begin{bmatrix} 1 & 0 \\ 0 & \eta \end{bmatrix} V \begin{bmatrix} 1 & 0 \\ 0 & \eta \end{bmatrix}^*$, it then follows that $V_{\eta} = \begin{bmatrix} a & B\eta^* \\ \eta C & \eta D \eta^* \end{bmatrix}$ is an isometry on $\bigoplus\limits_{i=1}^n (\clm_i \oplus \cln_i)$. Moreover, since $B \eta^* = \begin{bmatrix} B(1) & B(2) \end{bmatrix}$ and $\eta C = \begin{bmatrix} C(1) & C(2) \end{bmatrix}^t$, we see that
\[
V_{\eta} = \begin{bmatrix} a & B(1) & B(2) \\ C(1) & D(1) & \frac{1}{a} C(1) B(2) \\ C(2) & 0 & D(2) \end{bmatrix} \in \clb\Big(\mathbb{C} \oplus (\bigoplus_{i=1}^n \clm_i) \oplus (\bigoplus_{i=1}^n \cln_i)\Big),
\]
where
\begin{equation}\label{eq3 - BiCi}
B(i) = \begin{bmatrix} B_1(i) & B_2(i) \end{bmatrix} \quad \mbox{and} \quad C(i) = \begin{bmatrix} C_1(i) \\ C_2(i) \end{bmatrix},
\end{equation}
for all $i = 1, 2$. We have now arrived at the setting of the proof of Theorem \ref{th-factor 2} (more specifically, compare $V_{\eta}$ with $V$ in \eqref{eq-V33}). Following the constructions of $V_1$ and $V_2$ in the proof of Theorem \ref{th-factor 2}, we set
\begin{equation}\label{eq-3 V1 V2}
\begin{cases}
& V_1 = \begin{bmatrix} \alpha & B(1) \\ \frac{1}{\beta} C(1) & D(1) \end{bmatrix} \in \clb\Big(\mathbb{C} \oplus (\bigoplus\limits_{i=1}^n \clm_i)\Big)
\\
& \phantom{xxxx}
\\
& V_2 = \begin{bmatrix} \beta & \frac{1}{\alpha} B(2) \\ C(2) & D(2) \end{bmatrix} \in \clb\Big(\mathbb{C} \oplus (\bigoplus\limits_{i=1}^n \cln_i)\Big),
\end{cases}
\end{equation}
where
\[
|\beta|^2=|a|^2 + C(1)^* C(1)= 1 - C(2)^* C(2) \quad \text{and}\quad \alpha = \frac{a}{\beta}.
\]
Since $a \neq 0$, it follows that $\alpha$ (and $\beta$ too) is a non-zero scalars. One may now proceed, similarly as in the proof of Theorem \ref{th-factor 2}, to see that $V_1$ and $V_2$ are isometries. Then, applying Theorem \ref{th-MNH factor} to the pair of isometries $V_1$ and $V_2$, we get the canonical pair of isometries $\tilde{V}_1$ and $\tilde{V_2}$ such that $\tau_{\tilde{V}_1 \tilde{V}_2} = \tau_{V_1} \tau_{V_2}$.
On the other hand, it follows directly from the construction of $\tilde{V}_1$ and $\tilde{V}_2$ (see \eqref{eq-Def V}) that $V = \tilde{V}_1 \tilde{V}_2$ and consequently, $\tau_V = \tau_{\tilde{V}_1 \tilde{V}_2} = \tau_{V_1} \tau_{V_2}$. We have therefore proved the following counterpart of Theorem \ref{thm-sufficient V1V2} for isometric colligations satisfying property $\clf(n)$.

\begin{thm}\label{thm-3 between}
Let $V \in \clb(\mathbb{C} \oplus (\displaystyle \bigoplus_{i=1}^n (\clm_i \oplus\cln_i)))$ be an isometric colligation, and let $V$ satisfies property $\clf(n)$. If $\tau_V(0) \neq 0$ and $V$ admits the representation as in \eqref{eq3 - Vblock} with $B$, $C$ and $D$ as in \eqref{eq3 - BC} and \eqref{eq3 - D}, respectively, then
\[
V_1 = \begin{bmatrix} \alpha & B(1) \\ \frac{1}{\beta} C(1) & D(1) \end{bmatrix} \in \clb\Big(\mathbb{C} \oplus (\bigoplus_{i=1}^n \clm_i)\Big)
\quad \mbox{and} \quad
 V_2 = \begin{bmatrix} \beta & \frac{1}{\alpha} B(2) \\ C(2) & D(2) \end{bmatrix} \in \clb\Big(\mathbb{C} \oplus (\bigoplus_{i=1}^n \cln_i)\Big),
\]
are isometric colligations where $B(i), C(i)$ and $D(i)$ are as in \eqref{eq3 - D1D2} and \eqref{eq3 - BiCi} and $\alpha$ and $\beta$ are non-zero scalars and satisfies the following conditions
\[
|\beta|^2 = |a|^2 + C(1)^* C(1) \quad \mbox{and} \quad \alpha = \frac{a}{\beta}.
\]
Moreover, $\tau_V = \tau_{V_1} \tau_{V_2}$.
\end{thm}

This along with Theorem \ref{th-MNH factor} yields the following classification of Schur-Agler class factors of Schur-Agler class functions in $\cls\cla(\D^n)$, $n \geq 1$:

\begin{thm}\label{th-factor 3}
Suppose $\theta\in \mathcal{SA}(\mathbb{D}^n)$, and suppose that $\theta(0)\neq 0.$ Then $\theta = \phi \psi$ for some $\phi, \psi \in \mathcal{SA}(\mathbb{D}^n)$ if and only if $\theta = \tau_V$ for some isometric colligation $V$ satisfying property $\clf(n)$.
\end{thm}

Given $\theta = \tau_V$ for some isometric colligation $V$ satisfying property $\clf(n)$, as presented above, we now know that $\theta = \vp \psi$ for some $\phi, \psi \in \mathcal{SA}(\mathbb{D}^n)$. If $V$ admits the representation as in \eqref{eq3 - Vblock}, then it follows moreover from \eqref{eq-3 V1 V2} that
\begin{equation}\label{eq-3 phi psi}
\begin{cases}
& \phi(\z) = \alpha + \frac{1}{\beta} B(1) (I_{\clm_1^n} - E_{\clm_1^n}(\z) D(1))^{-1} E_{\clm_1^n}(\z) C(1)
\\
&\psi(\z) = \beta + \frac{1}{\alpha} B(2) (I_{\cln_1^n} - E_{\cln_1^n}(\z) D(2))^{-1} E_{\cln_1^n}(\z) C(2) \qquad (\z \in \D^n).
\end{cases}
\end{equation}
The assumption that $\theta(0) \neq 0$ in the proof of the sufficient part will be discussed in Section \ref{sbusect theta 0}. Also see Subsection \ref{FmFn} for a natural connection between $\clf_m(n)$ and $\clf(n)$, $1 \leq m < n$.

\section{Functions vanishing at the origin}\label{sbusect theta 0}

As pointed out in Remark \ref{rem-1}, factorizations of functions vanishing at the origin reveals more detailed structural properties of associated colligation matrices. To this end, in this section, we present a complete description of the connection between isometric colligations and Schur-Agler factors of Schur-Agler class functions vanishing at the origin. The case of one variable Schur functions will serve well to illustrate the notation scheme for functions in several variables that we adopt.

Suppose $\theta \in \cls(\D)$, $\theta(0) = 0$ and $\theta = \vp \psi$ for some $\phi$ and $\psi$ in $\cls(\D)$. The following two cases can arise:

\NI\textsf{Case (i)} $\phi(0)=0$ and $\psi(0)\neq 0$: Let $\vp = \tau_{V_1}$ and $\psi = \tau_{V_2}$, where $V_1=\begin{bmatrix}0&Q\\R&S\end{bmatrix}\in\mathcal{B}(\mathbb{C}\oplus\mathcal{H}_1)$ and $V_2=\begin{bmatrix}x&Y\\Z&W\end{bmatrix}\in\mathcal{B}(\mathbb{C}\oplus\mathcal{H}_2)$. Therefore $\tilde{V}_1 = \left[\begin{array}{@{}c|ccc@{}}
0&Q&0\\ \hline R&S&0 \\ 0&0&I
\end{array}\right]$ and $\tilde{V}_2 = \left[\begin{array}{@{}c|ccc@{}}
x&0&Y \\ \hline 0&I&0\\Z&0&W
\end{array}\right]$ are isometries in $\clb(\mathbb{C} \oplus (\clh_1 \oplus \clh_2))$. On defining $V := \tilde{V}_1 \tilde{V}_2$, we have the isometry
\begin{equation}\label{eq-new5 V case1}
V = \left[\begin{array}{@{}c|ccc@{}}
0&B_1&0\\\hline C_1&D_1&D_2\\ C_2&0&D_4
\end{array}\right] \in \clb(\mathbb{C} \oplus (\clh_1 \oplus \clh_2)),
\end{equation}
where
\[
\left[\begin{array}{@{}c|ccc@{}}
0&B_1&0\\\hline C_1&D_1&D_2\\ C_2&0&D_4
\end{array}\right]
=
\left[\begin{array}{@{}c|ccc@{}}
0&Q&0\\\hline xR&S&RY\\ Z&0&W
\end{array}\right].
\]
We then have $C_1 = xR$ and $D_2 = RY$, and consequently the condition $R^*R=1$ yields
\[
C_1C_1^*D_2=|x|^2RR^*D_2=|x|^2RR^*RY=|x|^2RY=|x|^2D_2 = C_1^*C_1D_2,
\]
as $C_1^*C_1 = |x|^2 (>0)$. Moreover, with $V$ as in \eqref{eq-new5 V case1}, we compute $\tau_V$ as:
\[
\begin{split}
\tau_V(z) & = \begin{bmatrix}B_1 & 0 \end{bmatrix} \Big(I - z \begin{bmatrix} D_1 & D_2 \\ 0 & D_4 \end{bmatrix} \Big)^{-1} z \begin{bmatrix} C_1 \\ C_2 \end{bmatrix}
\\
& = z \begin{bmatrix}B_1 & 0 \end{bmatrix} \begin{bmatrix} (I - z D_1)^{-1} & (I - zD_1)^{-1} z D_2 (I - zD_4)^{-1}\\ 0 & (I - zD_4)^{-1} \end{bmatrix} \begin{bmatrix} C_1 \\ C_2 \end{bmatrix}
\\
& = z \begin{bmatrix} B_1 (I - z D_1)^{-1} & z B_1 (I - zD_1)^{-1} D_2 (I - zD_4)^{-1}\end{bmatrix} \begin{bmatrix} C_1 \\ C_2 \end{bmatrix},
\end{split}
\]
and so
\begin{equation}\label{eq5-tauV}
\tau_V(z) = \Big(z B_1 (I - z D_1)^{-1} \Big) \Big(C_1 + z D_2 (I - zD_4)^{-1} C_2 \Big) \quad \quad (z \in \D).
\end{equation}
Substituting the values of $B_1, C_i$, and $D_j$, $i = 1, 2$ and $j = 2, 4$, we have
\[
\tau_V(z) = (z B_1 (I - z D_1)^{-1}) (C_1 + z D_2 (I - zD_4)^{-1} C_2) = (z Q (I - z D_1)^{-1}) (x R + z RY (I - z W)^{-1} Z),
\]
and hence $\tau_V(z) = (z Q (I - zS)^{-1} R) (x + zY (I - zW)^{-1}Z)$ for all $z \in \D$, which implies that $\theta = \tau_V$. Thus, we have collected together all the necessary properties of the isometric colligation $V$ as:
\begin{equation}\label{eq-new5 CD case1}
C_1C_1^*D_2 = C_1^*C_1D_2 \quad \mbox{and} \quad C_1^*C_1 >0.
\end{equation}
Conversely, suppose $V$ is an isometric colligation as in \eqref{eq-new5 V case1}, let $\theta = \tau_V$ and let $V$ satisfies the conditions in \eqref{eq-new5 CD case1}. Let $x$ be a non-zero scalar such that $|x|^2 = C_1^*C_1$. Define $V_1 \in \clb(\mathbb{C} \oplus \clh_1)$ and $V_2 \in \clb(\mathbb{C} \oplus \clh_2)$ by
\[
V_1 =\begin{bmatrix} 0&B_1\\\frac{1}{x}C_1&D_1 \end{bmatrix}
\quad\mbox{and}\quad
V_2 =\begin{bmatrix} x&\frac{1}{\bar{x}}C_1^*D_2\\C_2&D_4 \end{bmatrix}.
\]
Note that $|x|^2=1-C_2^*C_2=C_1^*C_1$. A simple computation then shows that $V_1$ and $V_2$ are isometric colligations. Now we compute
\[
\tau_{V}(z)=zB_1(1-zD_1)^{-1}C_1+z^2B_1(1-zD_1)^{-1}D_2(1-zD_4)^{-1}C_2,
\]
and
\[
\tau_{V_1}(z)\tau_{V_2}(z)=zB_1(1-zD_1)^{-1}C_1+z^2B_1(1-zD_1)^{-1}\{\frac{1}{|x|^2}C_1C_1^*D_2\}(1-zD_4)^{-1}C_2.
\]
Thus, $\tau_{V}=\tau_{V_1}\tau_{V_2}$ where $\tau_V(0) = \tau_{V_1}(0) = 0$ and $\tau_{V_2}(0) \neq 0$.

\vspace{0.2in}

\NI\textsf{Case (ii)} $\phi(0) = \psi(0) = 0$: Suppose $\phi = \tau_{V_1}$ and $\psi = \tau_{V_2}$, where $V_1=\begin{bmatrix}0&Q\\R&S\end{bmatrix}\in\mathcal{B}(\mathbb{C}\oplus\mathcal{H}_1)$ and $V_2=\begin{bmatrix}0&Y\\Z&W\end{bmatrix}\in\mathcal{B}(\mathbb{C}\oplus\mathcal{H}_2)$ are isometric colligations. We associate with $V_1$ and $V_2$ the isometric colligation
\[
V=\begin{bmatrix}
0&Q&0\\R&S&0\\0&0&I
\end{bmatrix}\begin{bmatrix}
0&0&Y\\0&I&0\\Z&0&W
\end{bmatrix}=\left[\begin{array}{@{}c|ccc@{}}
0&Q&0\\\hline 0&S&RY\\ Z&0&W
\end{array}\right],
\]
in $\clb(\mathbb{C} \oplus \clh_1 \oplus \clh_2)$ and set
\begin{equation}\label{eq5-case 2 v}
V = \left[\begin{array}{@{}c|ccc@{}}
0&B_1&0\\\hline 0&D_1&D_2\\ C_2&0&D_4
\end{array}\right].
\end{equation}
Then, in view of \eqref{eq5-tauV}, it follows that $\theta = \tau_{V}$. Also we pick the essential properties of the isometric colligation $V$ as
\begin{equation}\label{eq5-case 2 xdd}
X^*X=1, \quad X^*D_1 =0, \quad \mbox{and} \quad D_2=XY,
\end{equation}
where $X = R$. Note that the first two equalities follows from the fact that $V_1$ is an isometry.

\NI To prove the converse, suppose $V$ is an isometric colligation as in \eqref{eq5-case 2 v}, $\theta = \tau_V$, $X \in \clb(\mathbb{C}, \clh_2)$ is an isometry, $Y \in \clb(\clh_2, \mathbb{C})$ and the conditions in \eqref{eq5-case 2 xdd} hold. Since $V^* V = I$, we have
\[
\begin{bmatrix} C_2^* C_2 & 0 & C_2^* D_4 \\ 0 & B_1 B_1^* + D_1^* D_1 & D_1^* D_2 \\ D_4^* C_2 & D_2^* D_1 & D_2^* D_2 + D_4 D_4^* \end{bmatrix} = I_{\mathbb{C} \oplus \clh_1 \oplus \clh_2},
\]
and hence $V_1:= \begin{bmatrix} 0&B_1\\X&D_1 \end{bmatrix} \in \clb(\mathbb{C} \oplus \clh_1)$ is an isometric colligation. Since $D_2 = XY$, $D_2^* D_2 = Y^* Y$, and hence $D_2^* D_2 + D_4^* D_4 = I$ yields $Y^*Y + D_4^* D_4 = I$. Thus $V_2: = \begin{bmatrix}0&Y\\C_2&D_4 \end{bmatrix} \in \clb(\mathbb{C} \oplus \clh_2)$ is an isometric colligation. For all $z \in \D$, we have
\[
\begin{split}
\tau_{V_1}(z)\tau_{V_2}(z) & = z^2B_1(1-zD_1)^{-1}XY(1-zD_4)^{-1}C_2,
\end{split}
\]
and, on the other hand, in view of \eqref{eq5-tauV}, we have
\[
\tau_V(z) = z^2B_1(1-zD_1)^{-1} D_2(1-zD_4)^{-1}C_2.
\]
This and $XY = D_2$ implies that $\theta= \tau_{V}=\tau_{V_1}\tau_{V_2}$. Thus we have proved the following:

\begin{thm}\label{thm-new 5 case1}
Suppose $\theta\in \cls(\mathbb{D})$ and $\theta(0)=0$. Then:

(1) $\theta = \phi \psi$ for some $\phi, \psi \in \cls(\mathbb{D})$ and $\psi(0)\neq 0$ if and only if there exists an isometric colligation
\[
V = \left[\begin{array}{@{}c|ccc@{}}
0&B_1&0\\ \hline C_1&D_1&D_2\\ C_2&0&D_4
\end{array}\right] \in \clb(\mathbb{C} \oplus (\clh_1 \oplus \clh_2)),
\]
such that $C_1C_1^*D_2=C_1^*C_1D_2$, $C_1^*C_1>0$, and $\theta = \tau_{V}$.

(2) $\theta = \phi \psi$ for some $\phi, \psi \in \cls(\mathbb{D})$ and $\phi(0)= 0=\psi(0)$ if and only if there exists an isometric colligation
\[
V = \left[\begin{array}{@{}c|ccc@{}}
0&B_1&0\\ \hline 0 &D_1&D_2\\ C_2&0&D_4
\end{array}\right] \in \clb(\mathbb{C} \oplus (\clh_1 \oplus \clh_2)),
\]
such that $\theta = \tau_{V}$, $X^*D_1=0$, and $D_2=XY$ for some $Y\in \mathcal{B}(\mathcal{H}_2,\mathbb{C})$ and isometry $X\in \mathcal{B}(\mathbb{C},\mathcal{H}_2)$.
\end{thm}	

The general case of functions vanishing at the origin in several variables (in $\cls \cla(\D^n)$ or $\clm_1(H^2_n)$) can be studied using the technique developed in the proof of Theorem \ref{thm-new 5 case1}. In particular, similar arguments allow us to obtain also a similar classification of factorizations for functions in $\cls \cla(\D^n)$ vanishing at the origin. We only state the result in the setting of Section \ref{sec-factor Fn} and leave out the details to the reader.

\begin{thm}
Suppose $\theta \in  \mathcal{AS}(\mathbb{D}^n)$ and $\theta(0) = 0$. Then:

(1) $\theta = \phi \psi$ for some $\phi, \psi \in \cls \cla(\mathbb{D}^n)$ and $\psi(0)\neq 0$ if and only if there exist Hilbert spaces $\{\clh_i\}_{i=1}^n$, $\{\clm_i\}_{i=1}^n$ and $\{\cln_i\}_{i=1}^n$ and an isometric colligation
\[
V=\begin{bmatrix}0&B\\C&D\end{bmatrix}
= \left[\begin{array}{@{}c|ccc@{}}
0&  B_{1}  & \cdots & B_{n} \\\hline
C_{1} & D_{11} & \cdots & D_{1n} \\
\vdots & \vdots  & \ddots & \vdots \\
C_{n} & D_{n1} & \cdots  & D_{nn}
\end{array}\right]\in \mathcal{B}(\mathbb{C}\oplus(\bigoplus_{i=1}^n\mathcal{H}_i))
\]
such that $\theta = \tau_V$ and $\mathcal{H}_k = \mathcal{M}_k\oplus\mathcal{N}_k$, $k = 1, \ldots, n$, and representing $B_i$, $C_i$ an $D_{ij}$ as
\[
B_i=\left[B_i(1), B_i(2) \right]\in\mathcal{B}(\mathcal{M}_i\oplus\mathcal{N}_i,\mathbb{C}), \quad C_i=\begin{bmatrix}C_i(1)\\C_i(2)\end{bmatrix}\in \mathcal{B}(\mathbb{C}, \mathcal{M}_i \oplus \mathcal{N}_i),
\]
and $
D_{ij}=\begin{bmatrix} D_{ij}(1) & D_{ij}(12)\\ D_{ij}(21) & D_{ij}(2) \end{bmatrix} \in \mathcal{B}(\mathcal{M}_j \oplus \mathcal{N}_j, \mathcal{M}_i \oplus \mathcal{N}_i)$, one has $B_i(2) =$, $D_{ij}(21) = 0$, and
\[
C(1)C(1)^*D(12)=C(1)^*C(1)D(12)\quad \mbox{and} \quad  C(1)^*C(1)>0,
\]
where $i, j = 1, \ldots, n$, and $C(1)=\begin{bmatrix} C_1(1)\\\vdots\\C_n(1) \end{bmatrix}$ and $D(12)= \begin{bmatrix} D_{ij}(12)\end{bmatrix}_{i,j=1}^n.$

(2) $\theta = \phi \psi$ for some $\phi, \psi \in \cls \cla(\mathbb{D}^n)$ and $\phi(0)= 0=\psi(0)$ if and only if there exist Hilbert spaces $\{\clh_i\}_{i=1}^n$, $\{\clm_i\}_{i=1}^n$ and $\{\cln_i\}_{i=1}^n$, an isometry $X\in \mathcal{B}(\mathbb{C}, \bigoplus\limits_{i=1}^n \clm_i)$, a bounded linear operator $Y\in \mathcal{B}(\bigoplus\limits_{i=1}^n \cln_i, \mathbb{C})$ and an isometric colligation
\[
V=\begin{bmatrix}0&B\\C&D\end{bmatrix}
= \left[\begin{array}{@{}c|ccc@{}}
0&  B_{1}  & \cdots & B_{n} \\\hline
C_{1} & D_{11} & \cdots & D_{1n} \\
\vdots & \vdots  & \ddots & \vdots \\
C_{n} & D_{n1} & \cdots  & D_{nn}
\end{array}\right]\in \mathcal{B}\Big(\mathbb{C}\oplus(\bigoplus_{i=1}^n\mathcal{H}_i)\Big),
\]
such that $\theta = \tau_V$ and $\mathcal{H}_k = \mathcal{M}_k\oplus\mathcal{N}_k$, $k = 1, \ldots, n$, and representing $B_i$, $C_i$ an $D_{ij}$ as	
\[
B_i=\left[B_i(1), B_i(2) \right]\in\mathcal{B}\Big(\mathcal{M}_i\oplus\mathcal{N}_i,\mathbb{C}\Big), \quad C_i=\begin{bmatrix}C_i(1)\\C_i(2)\end{bmatrix}\in \mathcal{B}\Big(\mathbb{C},\mathcal{M}_i\oplus\mathcal{N}_i\Big),
\]
and $D_{ij}=\begin{bmatrix} D_{ij}(1) & D_{ij}(12)\\ D_{ij}(21) & D_{ij}(2) \end{bmatrix} \in \mathcal{B}(\mathcal{M}_j\oplus\mathcal{N}_j, \mathcal{M}_i \oplus \mathcal{N}_i)$,
one has $B_i(2) = 0$, $C_i(1) = 0$, and
\[
D_{ij}(21) = 0,  \quad
D(12)=XY \quad \mbox{and} \quad X^*D(1)=0,
\]
where
\[
D(1)=\left[D_{ij}(1) \right]_{i,j=1}^n \in \mathcal{B}\Big(\bigoplus\limits_{p=1}^n \mathcal{M}_p\Big),
\quad \mbox{and} \quad
\quad D(12)= \begin{bmatrix} D_{ij}(12)\end{bmatrix}_{i,j=1}^n \in \mathcal{B}\Big(\bigoplus\limits_{p=1}^n \mathcal{N}_p, \bigoplus\limits_{p=1}^n \mathcal{M}_p\Big).
\]
\end{thm}


\newsection{Examples and remarks}\label{sec last}

This section is devoted to some concrete examples, further results and general remarks concerning Schur functions.

\subsection{One variable factors} Our interest here is to analyze Schur-Agler class functions in $\cls\cla(\D^n)$ which can be factored as a product of $n$ Schur functions. More specifically, let $\vp \in \mathcal{SA}(\mathbb{D}^n)$ and let $\vp(0)\neq 0$. Suppose $\vp(\z) = \prod\limits_{i=1}^n \vp_i(z_i)$, $\z \in \D^n$, for some $\vp_i \in \cls(\D)$, $i = 1, \ldots, n$. Then there exist isometric colligations $V_i = \begin{bmatrix} a_1 & \hat{B}_i \\ \hat{C}_i & \hat{D}_i \end{bmatrix} \in \clb(\mathbb{C} \oplus \clh_i)$ such that $\vp_i = \tau_{V_i},$
for all $i = 1, \ldots, n$. Let $a = \prod_{i=1}^n a_i,$
and define
\[
\hat{V}_1 = \begin{bmatrix} a_1 & \hat{B}_1 & 0 \\ \hat{C}_1 & \hat{D}_1 & 0 \\ 0 & 0 & I \end{bmatrix},  \quad  \tilde{V}_n = \begin{bmatrix} a_n & 0 & \hat{B}_n \\ 0 & I & 0 \\ \hat{C}_n & 0 & \hat{D}_n \end{bmatrix}
\quad \mbox{and} \quad
\hat{V}_i = \begin{bmatrix} a_i &  0 & \hat{B}_i & 0 \\0 & I & 0 & 0 \\ \hat{C}_i & 0 & \hat{D}_i & 0 \\ 0 & 0 & 0 & I \end{bmatrix},
\]
in $\clb(\mathbb{C} \oplus \clh_1 \oplus \clh_2^n)$, $\clb(\mathbb{C} \oplus \clh_1^{n-1} \oplus \clh_n)$, and
 $\clb(\mathbb{C} \oplus \clh_1^{i-1} \oplus \clh_i \oplus \clh_{i+1}^n)$  respectively and for all $1 < i < n$. Then $V = \prod_{i=1}^n \hat{V}_i,$
is an isometry in $\clb(\mathbb{C} \oplus \clh_1^n)$. Moreover, it follows that
\begin{equation}\label{eq-V 1 variable}
V = \begin{bmatrix} a & B \\ C & D \end{bmatrix} = \left[\begin{array}{@{}c|ccc@{}}
a & B_{1}  & \cdots & B_{n} \\\hline
C_{1} & D_{11} & \cdots & D_{1n} \\
\vdots & \vdots  & \ddots & \vdots \\
C_{n} & D_{n1} & \cdots  & D_{nn}
\end{array}\right]
\end{equation}
where
\[
B_i = (\prod_{k=1}^{i-1} a_{k}) \hat{B}_i, \; C_i = (\prod_{k=i+1}^n a_{k}) \hat{C}_i,
\;\mbox{and}\;
D_{ij} =
\begin{cases} \hat{D}_{i} & \mbox{if}~ i = j
\\
0 & \mbox{if}~ i > j
\\
(a_{i+1} \cdots a_{j-1}) \hat{C}_i \hat{B}_j & \mbox{if}~ i < j. \end{cases}
\]
Hence $a D_{i,j} = C_{i} B_{j},$
for all $1\leqslant i<j\leqslant n$. Then by repeated application of Theorem \ref{th-1st factor}, we have $\vp = \tau_V$. The converse, as stated below, follows directly from repeated applications of Theorem \ref{thm-sufficient V1V2}. We have thus proved the following theorem.

\begin{thm}
Suppose $\theta\in \mathcal{SA}(\mathbb{D}^n)$ and $\theta(0)\neq 0$. Then $\theta(\z)= \displaystyle \prod_{i=1}^n \theta_i(z_i)$, $\z \in \D^n$ for some Schur functions $\{\theta_i\}_{i=1}^n \subseteq \cls(\D)$ if and only if $\theta = \tau_V$ for some isometric colligation $V = \begin{bmatrix} a & B \\ C & D \end{bmatrix} = \left[\begin{array}{@{}c|ccc@{}}
a & B_{1}  & \cdots & B_{n} \\\hline
C_{1} & D_{11} & \cdots & D_{1n} \\
\vdots & \vdots  & \ddots & \vdots \\
C_{n} & D_{n1} & \cdots  & D_{nn}
\end{array}\right]$ on $\mathbb{C} \oplus \Big(\bigoplus\limits_{i=1}^n \clh_i\Big)$ such that $D_{ij} =
\begin{cases} {D}_{i} & \mbox{if}~ i = j
\\
0 & \mbox{if}~ i > j
\\
\frac{1}{a} {C}_i {B}_j & \mbox{if}~ i < j \end{cases}$.
\end{thm}

\subsection{Examples}\label{subsection-examples}

Here we aim at applying our results to some concrete examples.

\NI \textsf{Example 1:} Let $\phi \in \cls(\D)$ and $\phi = \tau_{V_0}$ for some isometric colligation $V_0 = \begin{bmatrix} a & B \\ C & D \end{bmatrix} \in \clb(\mathbb{C}\oplus\mathcal{H})$. Now we consider $\psi(z) = z^m$, $z \in \D$ and $m \in \mathbb{N}$. One then shows that
\[	
V_m=\left[
\begin{array}{c|cccc}
0&1&0&\cdots &0\\  \hline
0&0&1&\cdots &0\\
\vdots&\vdots&\vdots&\ddots&\vdots\\
0&0&0&\cdots&1\\
1&0&0&\cdots&0
\end{array}\right]
\in \clb(\mathbb{C}\oplus\mathbb{C}^m),
\]
is an isometric colligation and $\psi = \tau_{V_m}$. Set $\theta = \phi \psi = \tau_{V_0} \tau_{V_m}$. Then by Theorem \ref{th-MNH factor} (or more specifically, by \eqref{eq-Def V}) it follows that $\tau_V(z) = z^m \vp(z)$, $z \in \D$, where $V \in \clb(\mathbb{C} \oplus \clh \oplus \mathbb{C}^m)$ is an isometric colligation with the following representation
\[	
V=\left[
\begin{array}{c|cc|ccccc}
0&B&a&0&0&\cdots &0\\
\hline
0&D&C&0&0&\cdots &0
\\
0&0&0&1&0&\cdots&0\\
\hline
0&0&0&0&1&\cdots&0
\\
\vdots& \vdots& \vdots& \vdots& \vdots&\ddots&\vdots
\\
0&0&0&0&0&\cdots&1\\
1&0&0&0&0&\cdots&0
\end{array}\right] \in \clb(\mathbb{C}\oplus(\mathcal{H}\oplus\mathbb{C})\oplus \mathbb{C}^{m-1}).
\]

\vspace{0.2in}

\NI \textsf{Example 2:} Our second example concerns Blaschke factors: If $\lambda \in \D$, then the \textit{Blaschke factor} $b_{\lambda} \in \mbox{Aut}({\D})$ is defined by
\[
b_{\lambda}(z) = \frac{z - \lambda}{1- \bar{\lambda} z} \quad \quad (z \in \D).
\]
Now observe that, for each $\lambda \in \D$, the matrix $V_{\lambda} =\begin{bmatrix}
-\lambda & \sqrt{1-|\lambda|^2}
\\
\sqrt{1-|\lambda|^2}&\bar{\lambda}
\end{bmatrix} \in \clb(\mathbb{C} \oplus \mathbb{C})$ is an isometric colligation and $b_{\lambda} = \tau_{V_{\lambda}}.$
Now, suppose $\alpha, \beta \in \D$ and $\theta(\z) = b_{\alpha}(z_1) b_{\beta}(z_2)$, $\z \in \D^2$. Then Theorem \ref{th-1st factor} implies that $\theta = \tau_V,$
where
\[
V = \begin{bmatrix}
\alpha\beta & \sqrt{1-|\alpha|^2} & -\alpha \sqrt{1-|\beta|^2}
\\
-\beta \sqrt{1-|\alpha|^2} & \bar{\alpha} & \sqrt{1-|\alpha|^2} \sqrt{1-|\beta|^2}
\\
\sqrt{1-|\beta|^2} & 0 & \bar{\beta}
\end{bmatrix},
\]
is an isometric colligation in $M_{3}(\mathbb{C})$.

\subsection{On $\clf_m(n)$ and $\clf(n)$}\label{FmFn}

Let $1 \leq m < n$. Suppose $V \in \clb(\mathbb{C} \oplus \clh_1^m \oplus \clh_{m+1}^n)$ satisfies property $\clf_m(n)$. On account of Theorem \ref{thm-sufficient V1V2}, we have
\[
\tau_V(\z) = \tau_{V_1}(z_1, \ldots, z_m) \tau_{V_2}(z_{m+1}, \ldots, z_n) \quad \quad (\z \in \D^n),
\]
for some isometric colligations $V_1 \in \clb(\mathbb{C} \oplus \clh_1^m)$ and $V_2 \in \clb(\mathbb{C} \oplus \clh_{m+1}^n)$. Note that $\tau_{V_1} \in \cls \cla(\D^m)$ and $\tau_{V_2} \in \cls \cla(\D^{n-m})$. The above factorization and Theorem \ref{th-factor 3} further implies that $\tau_V = \tau_{\tilde{V}}$ for some isometric colligation $\tilde{V} \in \clb(\mathbb{C}\oplus \Big(\bigoplus\limits_{i=1}^n (\clm_i \oplus \cln_i)\Big))$ satisfying property $\clf(n)$. It is then natural to ask to what extent one can recover $\tilde{V}$ from $V$. To determine the isometric colligation $\tilde{V}$, we proceed as follows: First, we let
\begin{equation}\label{eq-5 V}
V=\left[\begin{array}{@{}c|ccc@{}}
a&B_1&\cdots &B_n\\\hline
C_1 & D_{11} & \cdots & D_{1n}
\\
\vdots&\vdots&\ddots&\vdots\\
C_n & D_{n1} & \cdots & D_{nn}
\end{array}\right]\in\mathcal{B}(\mathbb{C}\oplus\mathcal{H}_1^n),
\end{equation}
where $D_{ij}=0$ for $i=m+1,\dots,n$ and $j=1,\ldots,m$;  $aD_{ij}=C_i B_j$ for $i=1,\ldots,m$ and $j=m+1,\dots,n$. Let $\cll$ be a Hilbert space. Set
\[
\clk_i =
\begin{cases}
\clh_i \oplus \cll & \mbox{if}~ 1 \leq i \leq m
\\
\cll \oplus \clh_i & \mbox{if}~ m+1 \leq i \leq n.
\end{cases}
\]
We now define
\[
Y_i =
\begin{cases}
\begin{bmatrix} B_i & 0 \end{bmatrix} & \mbox{if}~ 1 \leq i \leq m
\\
\\
\begin{bmatrix} 0 & B_i \end{bmatrix} & \mbox{if}~  m+1 \leq i \leq n,
\end{cases}
\quad \quad
Z_i =
\begin{cases}
\begin{bmatrix} C_i \\ 0 \end{bmatrix} & \mbox{if}~ 1 \leq i \leq m
\\
\\
\begin{bmatrix} 0 \\ C_i \end{bmatrix} & \mbox{if}~  m+1 \leq i \leq n,
\end{cases}
\]
and
\[
W_{ij} =
\begin{cases}
\begin{bmatrix} D_{ij} & 0 \\ 0 & \delta_{ij} I_{\cll} \end{bmatrix} & \mbox{if}~ 1 \leq i, j \leq m
\\
\\
\begin{bmatrix} \delta_{ij} I_{\cll} & 0 \\ 0 & D_{ij} \end{bmatrix}  & \mbox{if}~  m+1 \leq i, j \leq n,
\end{cases}
\]
and
\[
W_{ij} =
\begin{cases}
\begin{bmatrix} 0 & D_{ij} \\ 0 & 0 \end{bmatrix} & \mbox{if}~ 1 \leq i \leq m, \; m+1\leq j \leq n
\\
\\
\begin{bmatrix} 0 & 0 \\ D_{ij} & 0 \end{bmatrix}  & \mbox{if}~  m+1 \leq i \leq n, 1 \leq j \leq n.
\end{cases}
\]
Then, after some manipulations, it follows that the isometric colligation
\begin{equation}\label{eq-5 tildeV}
\tilde{V} := \left[\begin{array}{@{}c|ccc@{}}
a & Y_1 & \cdots & Y_n
\\
\hline
Z_1 & W_{11} & \cdots & W_{1n}
\\
\vdots&\vdots&\ddots&\vdots\\
Z_n & W_{n1} & \cdots & W_{nn}
\end{array}\right]\in\mathcal{B}(\mathbb{C}\oplus\mathcal{K}_1^n),
\end{equation}
satisfies property $\clf(n)$ and $\tau_V = \tau_{\tilde{V}}$. More specifically, we have proved the following:

\begin{thm}
Suppose $1 \leq m < n$ and let $V$ satisfies property $\clf_m(n)$. If the representation of $V$ is given by \eqref{eq-5 V}, then $\tau_{V} = \tau_{\tilde{V}}$, where $\tilde{V}$ is given by \eqref{eq-5 tildeV} and satisfies property $\mathcal{F}(n)$.
\end{thm}

\subsection{Factorizations of multipliers on the ball}\label{sub-ball}

Here we are interested in factorizations of multipliers of the Drury-Arveson space on the unit ball $\mathbb{B}^n$ in $\mathbb{C}^n$ \cite{BTV}. However (and curiously, if not surprisingly), the computations involved in representing multiplier factors of multipliers of the Drury-Arveson space seem relatively simpler than that of the Schur-Agler class functions on the polydisc. We omit details here and present only the final result.

Recall that the \textit{Drury-Arveson space}, denoted by $H^2_n$, is the Hilbert space of holomorphic functions on $\mathbb{B}^n$ corresponding to the reproducing kernel (cf. \cite{BTV})
\[
k(\z, \w) = (1 - \sum\limits_{i=1}^n z_i \bar{w}_i)^{-1} \quad \quad (\z, \w \in \ball^n).
\]
A complex-valued function $\vp$ on $\ball^n$ is said to be a \textit{multiplier} if $\vp H^2_n \subseteq H^2_n$. If $\vp$ is a multiplier, then $M_{\vp} f \mapsto \vp f$, $f \in H^2_n$, defines a bounded operator on $H^2_n$. We let $\clm(H^2_n)$ denote the commutative Banach algebra of multipliers equipped with the operator norm $\|\vp\| := \|M_{\vp}\|_{\clb(H^2_n)}$. Also we define
\[
\clm_1(H^2_n) = \{\vp \in \clm(H^2_n): \|{\vp}\| \leq 1\}.
\]
The following characterization of multipliers (see \cite{BTV, EP}), parallel to the transfer function realizations of Schur-Agler class functions on $\D^n$ (see Theorem \ref{thm-Agler Dn}), is the starting point: Suppose $\vp$ is a complex-valued function on $\ball^n$. Then $\vp \in \clm_1(H^2_n)$ if and only if there exist a Hilbert space $\clh$ and an isometric colligation $V = \begin{bmatrix} a & B \\ C & D \end{bmatrix} : \mathbb{C} \oplus \clh \raro \mathbb{C} \oplus \clh^n$ such that $\vp = \tau_V$, where
\[
\tau_V(\z) = a + B (I_{\clh} - E_{\clh^n}(\z) D)^{-1} E_{\clh^n}(\z) C \quad \quad (\z \in \mathbb{B}^n).
\]
Here given a Hilbert space $\clh$, we denote by $\clh^n$ the $n$-copies of $\clh$, and $E_{\clh^n} : \ball^n \raro \clb(\clh^n, \clh)$ the row operator $E_{\clh^n}(\z) = (z_1 I_{\clh}, \ldots, z_n I_{\clh})$, $\z \in \ball^n$.

We omit the proof of the following result which is similar (in spirit) to the case of $\cls \cla(\D^n)$.

\begin{thm}
Suppose $\theta \in \mathcal{M}_1(H_n^2)$ and $\theta(0)\neq 0$. There exist multipliers $\phi\in \mathcal{M}_1(H_m^2)$ and $\psi\in \mathcal{M}_1(H_{n-m}^2)$ such that $\theta(\z)=\phi(z_1,\ldots,z_m) \psi(z_{m+1}, \ldots, z_n)$, $\z \in \D^n$, if and only if  $\theta = \tau_V$ for some isometric colligation
\[
V = \begin{bmatrix} a & B \\ C & D \end{bmatrix}: \mathbb{C} \oplus (\clh_1 \oplus \clh_2) \raro \mathbb{C} \oplus (\clh_1 \oplus \clh_2)^n,
\]
such that writing $B = \begin{bmatrix} B(1) & B(2) \end{bmatrix}$, $C = \begin{bmatrix} C_1 & \ldots & C_n \end{bmatrix}^t$ and $D = \begin{bmatrix} D_1 & \ldots & D_n \end{bmatrix}^t$, one has
\[
C_j =
\begin{cases}
\begin{bmatrix} C_j(1) \\ 0 \end{bmatrix} & \mbox{if~} 1 \leq j \leq m
\\
\\
\begin{bmatrix} 0 \\ C_j(2) \end{bmatrix} & \mbox{if~} m+1 \leq j \leq n, \end{cases}
\qquad
D_j =
\begin{cases}
\begin{bmatrix} D_j(1) & D_j(2) \\ 0 & 0 \end{bmatrix} & \mbox{if~} 1 \leq j \leq m
\\
\\
\begin{bmatrix} 0 & 0 \\ 0 & D_{j}(3) \end{bmatrix} & \mbox{if~} m+1 \leq j \leq n, \end{cases}
\]
and $a D_i(2)=C_i(1)B(2)$ for all $i = 1, \ldots, m$.
\end{thm}

\subsection{Reversibility of factorizations}\label{sub-reverse}

A natural question to ask in connection with Theorem \ref{th-factor 3} is whether the canonical constructions of the colligation $V$ (out of a pair of isometric colligations $V_1$ and $V_2$) satisfying property $\clf(n)$ as in \eqref{eq-Def V} and $V_1$ and $V_2$ (out of an isometric colligation $V$ satisfying property $\clf(n)$) as in \eqref{eq-3 V1 V2} are reversible.

To answer this, we proceed as follows: Given $n \in \mathbb{N}$, we let $C(n)$ denote the set of all isometric colligations of the form $\begin{bmatrix} a & B \\ C & D \end{bmatrix} \in \clb(\mathbb{C} \oplus \Big(\bigoplus\limits_{i=1}^n \clh_i\Big))$ for some Hilbert spaces $\{\clh_i\}_{i=1}^n$, and let $F(n)$ denote the set of all isometric colligations satisfying property $\clf(n)$. Define $\pi : C(n) \times C(n) \raro F(n)$ by \[
\pi(V_1, V_2) = V \quad \quad (V_1, V_2 \in C(n)),
\]
where $V$ is as in \eqref{eq-Def V} (or Theorem \ref{th-MNH factor}). Also define $\kappa : F(n) \raro C(n) \times C(n)$ by
\[
\kappa(V) = (V_1, V_2) \quad \quad (V \in F(n)),
\]
where $V_1$ and $V_2$ are as in \eqref{eq-3 V1 V2}. Given $V_1$ and $V_2$ in $C(n)$, the aim here is to compare $\kappa (\pi(V_1, V_2))$ with $(V_1, V_2)$. Suppose $V_1 = \begin{bmatrix} \alpha & B \\ C & D \end{bmatrix} \in \clb \Big(\mathbb{C} \oplus (\displaystyle \bigoplus_{i=1}^n \clm_i)\Big)$ and $V_2 = \begin{bmatrix} \beta & F \\ G & H \end{bmatrix} \in \clb\Big(\mathbb{C} \oplus  (\displaystyle \bigoplus_{i=1}^n \cln_i)\Big)$ are isometric colligations and $a = \alpha \beta \neq0$. Then by \eqref{eq-Def V}, it follows that
\[
\pi(V_1, V_2) \in \clb\Big(\mathbb{C} \oplus (\displaystyle \bigoplus_{i=1}^n (\clm_i \oplus \cln_i))\Big),
\quad\mbox{and}\quad
\pi(V_1, V_2) = \left[\begin{array}{@{}cccc@{}}
\alpha \beta & \hat{B}_{1}  & \cdots & \hat{B}_{n} \\
\hat{C}_{1} & \hat{D}_{11} & \cdots & \hat{D}_{1n} \\
\vdots & \vdots  & \ddots & \vdots \\
\hat{C}_{n} & \hat{D}_{n1} & \cdots  & \hat{D}_{nn}
\end{array}\right],
\]
where $\hat{B}_i$, $\hat{C}_i$ and $\hat{D}_{ij}$, $i, j = 1, \ldots, n$, are given by as in \eqref{eq-3.1 one}. Since $\pi(V_1, V_2)$ satisfies property $\clf(n)$, in view of \eqref{eq-3 V1 V2}, it follows that $\kappa (\pi(V_1, V_2)) = (\tilde{V}_1, \tilde{V}_2)$, where
\[
\tilde{V}_1 = \left[\begin{array}{@{}cccc@{}}
\tilde{\alpha} & {B}_{1}  & \cdots & {B}_{n} \\
\frac{\beta}{\tilde{\beta}} {C}_{1} & {D}_{11} & \cdots & {D}_{1n} \\
\vdots & \vdots  & \ddots & \vdots \\
\frac{\beta}{\tilde{\beta}} {C}_{n} & {D}_{n1} & \cdots & {D}_{nn}
\end{array}\right] \quad \mbox{and} \quad
\tilde{V}_2 = \left[\begin{array}{@{}cccc@{}}
\tilde{\beta} & \frac{\alpha}{\tilde{\alpha}} F_{1}  & \cdots & \frac{\alpha}{\tilde{\alpha}} F_{n} \\
G_{1} & H_{11} & \cdots & H_{1n} \\
\vdots & \vdots  & \ddots & \vdots \\
G_{n} & H_{n1} & \cdots & H_{nn}
\end{array}\right],
\]
and $\tilde{\alpha}$ and $\tilde{\beta}$ are non-zero scalars satisfying the following relations
\[
|\tilde{\beta}|^2 = |\alpha|^2 |\beta|^2 + |\beta|^2 \Big(\sum_{i=1}^n C_i^* C_i\Big) \quad \mbox{and} \quad \tilde{\alpha} = \frac{\alpha \beta}{\tilde{\beta}}.
\]
But we know from $V_1^* V_1 = I$ that $|\alpha|^2 + C^* C = 1$, that is $|\alpha|^2 + \sum_{i=1}^n C_i^* C_i = 1.$
So $\tilde{\beta} =\bar{\varepsilon} \beta$ and $\tilde{\alpha} = \varepsilon \alpha$ for some unimodular constant $\varepsilon$. Hence
\[
\kappa \circ \pi \Big(\begin{bmatrix} \alpha & B \\ C & D \end{bmatrix}, \begin{bmatrix} \beta & F \\ G & H \end{bmatrix} \Big) =  \Big(\begin{bmatrix} \varepsilon \alpha & B \\ \varepsilon C & D \end{bmatrix}, \begin{bmatrix} \bar{\varepsilon} \beta & \bar{\varepsilon} F \\ G & H \end{bmatrix} \Big),
\]
where $\varepsilon$ is an unimodular constant.

One could equally consider the same question for Theorem \ref{th-factor 2}. The answer is similar and we leave the details to the reader.




\vspace{0.2in}

\noindent Data Availability: All data generated or analysed during this study are included in this published article.

\smallskip
\smallskip

\noindent Acknowledgement: The second author is supported in part by NBHM (NBHM/R.P.64/2014), and the Mathematical Research Impact Centric Support, MATRICS (MTR/2017/000522), and Core Research Grant (CRG/2019/000908), by the Science and Engineering Research Board (SERB), Department of Science \& Technology (DST), Government of India.

\end{document}